\documentclass[12pt,a4paper]{amsart}

\usepackage{amsmath,amsthm,amsfonts,amssymb,mathrsfs}
\usepackage{subfigure} 
\usepackage{graphicx}
\usepackage{color}

\usepackage{array}
\usepackage{booktabs}
\usepackage{caption}
\usepackage{float}

\usepackage{amsmath}
\usepackage{amssymb}
\usepackage{epsfig}
\usepackage{float}
\usepackage{latexsym,graphics}
\usepackage{dsfont}
\newtheorem{theorem}{Theorem}[section]
\newtheorem{proposition}[theorem]{Proposition}
\newtheorem{lemma}[theorem]{Lemma}

\newtheorem{definition}[theorem]{Definition}

\newtheorem{remark}[theorem]{Remark}

\newcommand{\dst}{\displaystyle}

\newcommand{\R}{\mathbb{R}}

\def\sqr#1#2{{\vcenter{\vbox{\hrule height.#2pt\hbox{\vrule width.#2pt height#1pt \kern#1pt\vrule width.#2pt}\hrule height.#2pt}}}}

\newcommand{\modulo}[1]{{\left|#1\right|}} 	
\newcommand{\norma}[1]{{\left \|#1\right \|}} 		
\renewcommand{\phi}{\varphi}		

\newenvironment{proofof}[1]
{\smallskip\noindent{\textbf{Proof~of~#1.}}
\hspace{1pt}}{\hspace{-5pt}{\nobreak\nobreak\hfill\nobreak
$\square$\vspace{2pt}\par}\smallskip\goodbreak}

\allowdisplaybreaks[4]

\begin{document}

\title[Analysis of particle methods for structured population models]{Analysis of particle methods for structured population models with nonlocal boundary term in the framework of bounded Lipschitz distance}

\author[P.Gwiazda]{Piotr Gwiazda$^1$}
\author[J.Jablonski]{Jedrzej Jablonski$^1$}
\author[A.Marciniak-Czochra]{Anna Marciniak-Czochra$^2$}
\author[A.Ulikowska]{Agnieszka Ulikowska$^1$}
\address[$^1$]{
Institute of Applied Mathematics and Mechanics,  University of Warsaw, Warszawa 02-097, Poland}
\email{pgwiazda@mimuw.edu.pl}
\address[$^2$]{
Institute of Applied Mathematics,  Interdisciplinary Center for Scientific Computing (IWR) and BIOQUANT, University of Heidelberg, 69120 Heidelberg, Germany}
\email{anna.marciniak@iwr.uni-heidelberg.de}

\date{\today}

\begin{abstract}
Recently developed theoretical framework for analysis of structured population dynamics in the spaces of nonnegative Radon measures with a suitable metric provides a rigorous tool to study numerical schemes based on particle methods. The approach is based on the idea of tracing growth and transport of measures which approximate the solution of original partial differential equation. In this paper we present analytical and numerical study of two versions of Escalator Boxcar Train (EBT) algorithm which has been widely applied in theoretical biology, and compare it to the recently developed split-up algorithm. The novelty of this paper is in showing well-posedness and convergence rates of the schemes using the concept of semiflows on metric spaces. Theoretical results are validated by numerical simulations of test cases, in which distances  between simulated and exact solutions are computed using flat metric.\\


\end{abstract}
\maketitle

\section{Introduction}

Nonlinear first order hyperbolic equations are often used in applications to describe crowd dynamics, for example pedestrian flows \cite{Bellomo,Bellomo1,Bellomo2,Piccoli1}, or dynamics of  populations which are  heterogenous in respect to some individual property (physiologically structured population models) \cite{diekmanngetto,diekmannmetz,DMCPZ,GJMC,webb}. The problems often lead to the models with nonlocal terms reflecting the impact of the whole population on birth and death processes of each individual \cite{diekmannmetz,GMC,perthame,webb} or on the velocity of the individual movement \cite{Muntean,Piccoli2,Piccoli1}. 

Recently, such problems have been also studied using methods originated from kinetic theory  and related to them metrics, such as Wasserstein and Monge-Kantorovich metrics \cite{Ackleh1,Ackleh2,CCGU,Muntean,GLMC,GMC, Piccoli2,two_sex}. The idea of representing a heterogeneous population as a sum of masses concentrated in different points of the structure follows from empirical understanding of the data obtained in measurements taken place in discrete time points. The choice of spaces of positive Radon measures  was proposed in  \cite{diekmannmetz} as relevant for biological applications, in which the initial distribution of individuals is concentrated with respect to the structure, i.e. is not absolutely continuous with respect to the Lebesgue measure.

 A framework for the analysis of solutions of structured population models using Wasserstein type metrics, adjusted to the nonconservative character of the considered problem, has been proposed in \cite{GLMC} using a flat metric (bounded Lipschitz distance) and in \cite{GMC} using a Wasserstein type metric adjusted to spaces of nonnegative Radon measures with integrable first moment. The advantage of that approach is in providing a structure of a space appropriate to compare solutions and to study their stability. Among others, continuous dependence with respect to the modeling ingredients is important in the context of numerical approximation and model calibration based on experimental data. 

In the current paper we focus on the classical  nonlinear renewal equation with a nonlocal boundary condition given by
\begin{eqnarray}\label{eq1}
\partial_t \mu_t + \partial_{x}(b(t,\mu_t)\mu_t) + c(t,\mu_t) \mu_t
& = &
0
\\ \nonumber
b(t,\mu_t)(x_b)D_{\lambda}\mu_t(x_b^+)
& = &
\int_{x_b}^{+\infty} \beta(t,\mu_t)(x) \; d \mu_t(x)
\end{eqnarray}
with nonnegative parameters $c$ and $\beta$ and a strictly positive $b$.

This model has served as a test case for the analysis of solutions in spaces of positive Radon measures developed in \cite{GLMC} and \cite{GMC}. The result on Lipschitz continuous dependence of solutions on the model parameters and initial data was then extended using a split-up method to the size-structured population models, i.e involving nonlocal terms in the equation \cite{CCGU}, and to systems of equations in \cite{two_sex}. As mentioned in \cite{GMC} the constructive proofs of Lipschitz continuous dependence of measure-valued solutions in respect to time and initial data provide a good tool to study convergence of the numerical algorithms based on particle methods.  

Recently, a numerical approach based on the split-up algorithm was proposed in \cite{CGU}. The split-up algorithm is based on the idea  of  representing a semigroup defining the model solution as a product of two semigroups related to equations with simpler structure. In case of model \eqref{eq1} it allows for reduction of the  problem involving transport terms and nonlocal growth terms to two problems involving either one or the other process. Analytical framework obtained in \cite{CCGU} allows to control convergence rate of the algorithm.

Another numerical approach based on measure-valued solutions of model \eqref{eq1} was proposed in \cite{EBT} using so called Escalator Boxcar Train (EBT) algorithm. The method is based on representing the solution as a sum of masses localized in discrete points and tracing their spatio-temporal evolution. 
The algorithm has been used in applications since a long time \cite{Roos}, however convergence of the scheme was shown only recently \cite{EBT} using the approach of metric spaces proposed in \cite{GLMC}. 

Similar approaches based on particle methods have been applied to simulate kinetic models from physics since more than three decades, see \cite{Harlow,Issautier, Raviart,Tadmor} and references therein. Recent examples of applications include porous medium equation \cite{Westdickenberg3},  isentropic Euler equations \cite{Westdickenberg1,Westdickenberg2}, and pedestrian flows \cite{Muntean}.  

Challenges arising from applications in biology using structured population models of type \eqref{eq1} are related to the nonlocal character of the boundary condition, which induces a regularization effect leading to solutions which are absolutely continuous with respect to the Lebesgue measure. Such solutions  need to be approximated by a sum of Dirac measures to allow further application of the particle method. 

In this paper, we provide a systematic numerical approach to problem \eqref{eq1} based on particle methods. 
In Section \ref{ParticleMethods} we present numerical schemes.  They include two versions of  EBT algorithm with different implementation of the boundary condition, as proposed originally in \cite{Roos} and then simplified in \cite{EBT}. Additionally, we present split-up algorithm, as proposed in \cite{CGU}. Our motivation is to compare different approaches.
The novelty of this work is in analysis of the two versions of EBT method. Since the original EBT algorithm is based on solving a system of ordinary differential equations with non-Lipschitz and even non-continuous right hand side, its well-posedness does not follow from a standard theory. We cope with this difficulty and provide a proof of local existence of solutions of the method. According to our knowledge this problem has not been solved before. Furthermore, we show an example in which the local solutions cannot be extended in time and a singularity appears.

Importantly, our analytical approach allows to control convergence rate of the algorithms, which is presented in Section \ref{Convergence}. The new proofs are based on approach of semiflows on metric spaces developed in \cite{ColomboGuerra2009} and applied to structured population models in \cite{GMC}. It allows for calculation of the estimates controlling the rate of convergence of the algorithms, which has not been obtained in previous studies of the EBT method. To validate practical applicability of the analytical results, we have implemented the EBT schemes and, also, the split-up algorithm. In Section \ref{Simulations} we present numerical simulations of test cases taken from literature, for which exact solutions can be obtained analytically. This allows us to check the accuracy of the numerical methods. The solutions are compared using flat distance, which is computed by a recently developed algorithm based on linear programming \cite{JMC}.  Comparing the results of the three numerical approaches we observe that the most efficient is the original EBT method, although the differences between the methods are not too large.

\section{Notation}\label{Notation}
Let $\R_{\geq} = \{x \in \R  :  0 \leq x\}$,
$\mathcal M^+(\R_{\geq})$ be the space of finite, nonnegative Radon measures and
$C_{b}^{\alpha,1} \left(
 [0,T]\times \mathcal M^+(\R_{\geq}) ;\; X \right)$ denote the space of $X$ valued functions, bounded with respect to the $\|\cdot\|_X$ norm, H\"older continuous with respect to time and Lipschitz continuous with respect to the second variable.

For each $f \in C_{b}^{\alpha,1} \left(
 [0,T]\times\mathcal M^+(\R_{\geq}) ;\; X \right)$ we define
$$
\norma{f}_{BC} = \sup_{(t,\mu) \in [0,T]\times \mathcal M^+(\R_{\geq}) } \norma{f(t,\mu)}_{X}
+
\sup_{t \in [0,T]} \mathrm{Lip}(f(t,\cdot))
+
\sup_{\mu \in \mathcal M^+(\R_{\geq}) } H(f(\cdot, \mu)),
$$
where $\mathrm{Lip}(f)$ denotes the Lipschitz constant and $H(f)$ is equal to
$$
H(f) = \sup_{s_1,s_2 \in [0,T]} \frac{\modulo{f(s_1) - f(s_2)}}{\modulo{s_1 - s_2}^{\alpha}}.
$$
The space $W^{1,\infty}$ is equipped with its usual norm, i.e.
$$
\norma{\gamma}_{W^{1,\infty}} = \max\{ \norma{\gamma}_{L^{\infty}}, \norma{\partial_x \gamma}_{L^{\infty}}\}.
$$
We define flat metric $\rho_F$ (also known as
bounded Lipchitz distance) as a distance derived from the dual norm of
$W^{1,\infty}$ (see e.g. \cite{MullerOrtiz}, \cite{Zhidkov}).
\begin{definition}\label{defflatmetric} 
Let $\mu, \nu \in{\mathcal M^+}({\R_{\geq}})$.  The distance
function $\rho_F : {\mathcal M^+}({\R_{\geq}}) \times  {\mathcal M^+}({\R_{\geq}}) \rightarrow
[0, \infty]$ is defined by
\begin{equation}\label{Flatmetric}
  \rho_F(\mu, \nu) :=
\dst \sup \Big\{\int_{\R_{\geq}} \psi  d (\mu - \nu)
  \big| \: \psi \in C^1({\R_{\geq}}), \|\psi\|_{W^{1,\infty}} \leq 1 \Big\},
  \end{equation}
where $\|\psi\|_{W^{1,\infty}} = \max\{\|\psi\|_{\infty},\|\partial_x \psi\|_{\infty}\}.$
\end{definition}
\noindent

\begin{definition}\label{DefSemigroup}
Let $(E,\rho)$ be a metric space. A Lipchitz semiflow is a semigroup $S: [0, \delta] \times [0,T] \times E \rightarrow E$
satisfying $$\rho\big(S(t; \tau)\mu,S(s;\tau)\nu\big) \leq L\left( \rho(\mu,\nu)
 + |t - s|\right),$$
where $ s,t\in [0, \delta]$, $\tau, s+\tau, t+\tau \in [0, T]$ and $\mu, \nu \in E$.
\end{definition}
\noindent

\section{Presentation of the algorithms}\label{ParticleMethods}

A concept of the particle methods is based on tracking groups of individuals which are similar to each other. Following the commonly applied terminology in natural sciences we call these groups cohorts. They are characterized by masses and locations $\{(m^i_o, x^i_o)\}_{i=B}^N$. Formally, the initial distribution $\mu_o$ is approximated by a sum of Dirac measures $\mu_o^n := \sum_{i= B}^N m^i_o \delta_{x^i_o}$, where $n = N - B + 1$.  Once the individual is allocated in the particular cohort, it stays there till the moment of death. Properties of each cohort change in time.  $m^i(\cdot)$ changes its value due to the growth processes, while $x^i(\cdot)$ evolves along the characteristics defined by the transport operator. New cohorts, which account for the influx of new individuals due to the birth processes, appear at the boundary with a fixed frequency, usually once per a time step. A cohort which currently corresponds to this phenomenon is called the boundary cohort.
Output of the algorithms at time $t$ is given as a sum of Dirac deltas $\mu_t^n := \sum_{i =B}^{N+K}{m^i(t)\delta_{x_i(t)}}$, where $K$ is the number of boundary cohorts created till the time moment $t$ and $n$ is the parameter of the initial data approximation. If no ambiguity occurs, we omit the superscript $n$ in the present section. To represent the result as a density function, it is necessary to mollify the output measure, for example, by using the formula given by equation \eqref{density}.

\noindent

\subsection{Escalator Boxcar Train (EBT)}\label{subsection_EBT_org}

The EBT algorithm allows to compute $\mu_t$, i.e. the approximate solution to \eqref{eq1} at time $t$ for the initial data $\mu_0$. The procedure iteratively solves a system of ODEs on a sufficiently short time interval $[ t_k, t_{k+1}]$. The output measure  $\mu_{ t_{k+1}}$ provides an input to  the next step of the algorithm.
Each step of the algorithm is based on introducing a new empty cohort and then  solving the following ODE system 
\begin{eqnarray}\label{EBT_org_system}
\left\{
	\begin{array}{rcll}
\displaystyle
\mathrm{\frac{d}{dt}} x^i(t) &=& b(t, \mu)(x^i(t)),&
\mathrm{for}\;\; i=B+1, \dots, J,
\\[2mm]
\displaystyle
\mathrm{\frac{d}{dt}} m^i(t) &=& -c(t, \mu)(x^i(t))m^i(t),
&
\mathrm{for}\;\; i=B+1, \dots, J,
	\end{array}
\right.
\end{eqnarray}

\begin{eqnarray}\label{EBT_org_BC}
\left\{
	\begin{array}{rcl}
\displaystyle
x^B(t) &=& 
\left\{
\begin{array}{lcl}
\frac{\pi^B(t)}{m^B(t)} + x_b,
&&
\mathrm{if}\;\; m^B(t) > 0,
\\
x_b, && \mathrm{otherwise}.
\end{array}
\right.\\ \\
\mathrm{\frac{d}{dt}} \pi^B(t) &=& b(t, \mu)(x_b)m^B(t)
+ \partial_x b( t,\mu )(x_b)\pi^B(t)
\\
&& - c( t, \mu)(x_b)\pi^B(t),
\\[3mm]
\displaystyle
\mathrm{\frac{d}{dt}} m^B(t) &=&-c(t, \mu)(x_b)m^B(t)  - \partial_x c (t, \mu)(x_b)\pi^B(t) 
\\
&& + \sum_{i=B}^{J} \beta(t, \mu)(x^i(t))m^i(t),\\
m^B (t_k) &=& 0,\\
\pi^B (t_k) &=& 0.
	\end{array}
\right.
\end{eqnarray}

Index $B$ indicates the boundary cohort. The other, so called internal, cohorts are denoted by $i = B+1, \dots, J$. The solution is given by $\mu_t = \sum_{i=B}^{J} m^i(t) \delta_{x^i(t)}$.

We assume that
\begin{equation}\label{assumptions:sim}
b,  c,  \beta \in
C_{b}^{\alpha,1} \left(
[0,T]\times \mathcal M^+(\R_{\geq} );
W^{1,\infty}(\R_{\geq} ) 
\right).
\end{equation}
%

\begin{theorem}\label{existence_org}
Assume that \eqref{assumptions:sim}  holds and the initial data in \eqref{EBT_org_system}~-~\eqref{EBT_org_BC} are nonnegative. Then, there exists a local in time, unique solution $(x^i,m^i)$, which is nonnegative and Lipschitz continuous in respect to time.
\end{theorem}

%
%
The EBT algorithm is derived using a linearisation of equations for total mass $m= \int_{\Omega_B(t)} u(t,x) dx$ and centre of  mass $\pi= \int_{\Omega_B(t)}x u(t,x) dx $. Because $\int_{\Omega_B(t)}x^2 u(t,x) dx$ cannot be expressed in terms of $m$ and $\pi$, the whole product $x c(t,\mu)(x)$ is linearised (for details see the derivation of the algorithm in Appendix). This truncation unexpectedly leads to lack of global in time existence of solutions and possible blow-up of $x^B$, as shown in Table \ref{Table5}. 
\begin{remark}
One can consider a modification of the EBT algorithm (EBT$^*$), which differs by the linearisation procedure of  $\int_{\Omega_B(t)}x c(t,x) u(t,x) dx$.
Linearising $c(t,\mu)(x)$ and taking $\int_{\Omega_B(t)}x^2u(t,x) dx =\int_{\Omega_B(t)}x^2 m^B(t) d \delta_{x^B(t)}(x)$, we obtain a version of the EBT algorithm, where equations on $\pi^b$ and $m^B$ in \eqref{EBT_org_BC} are replaced by
\begin{eqnarray}\label{EBT*}
\left\{
	\begin{array}{rcl}
\displaystyle
\mathrm{\frac{d}{dt}} \pi^B(t) &\!\!\!=\!\!\!& b(t, \mu)(x_b)m^B(t) - \partial_x c (t, \mu)\pi^B(t)x^B(t) 
+ \partial_x b( t,\mu )(x_b)\pi^B(t)
\\
&& - c( t, \mu)(x_b)\pi^B(t),
\\[3mm]
\displaystyle
\mathrm{\frac{d}{dt}} m^B(t) &\!\!\!=\!\!\!&-c(t, \mu)(x_b)m^B(t)  - \partial_x c (t, \mu)(x_b)\pi^B(t) 
\\
&& + \sum_{i=B}^{J} \beta(t, \mu)(x^i(t))m^i(t).
	\end{array}
\right.
\end{eqnarray}
This correction may lead to better properties of the algorithm including nonegativity of mass and a control of the boundary cohort centre $x^B<x^1$ (see Subsection \ref{tc3} and Table \ref{Table5}). However, in this paper we do not study this ODEs system rigorously.
\end{remark}

\subsection{Escalator Boxcar Train with simplified boundary equations (sEBT)}\label{subsection_sEBT}

The EBT algorithm can be modified, as proposed in \cite{EBT}, by replacing equation \eqref{EBT_org_BC} by
\begin{eqnarray}\label{EBT_simply_BC}
\left\{
	\begin{array}{rcl}
\displaystyle
\mathrm{\frac{d}{dt}} x^B(t) &=& b(t, \mu)(x^B(t)),
\\[2mm]
\displaystyle
\mathrm{\frac{d}{dt}} m^B(t) &=& -c( t, \mu)(x^B(t))m^B(t) + \sum_{i=B}^{J} \beta( t, \mu)(x^i(t))m^i(t).
	\end{array}
\right.
\end{eqnarray}

\begin{theorem}\label{existence_org2}
Assume that \eqref{assumptions:sim} holds and the initial data in  \eqref{EBT_org_system}, \eqref{EBT_simply_BC} are nonnegative. Then, there exists a global in time, unique solution $(x^i,m^i)$, which is nonnegative and Lipschitz continuous in respect to time.
\end{theorem}

\subsection{Split up algorithm (SU)}

\noindent
Split-up algorithm, proposed recently in \cite{CCGU},  iteratively provides approximations $\mu_{t_k}$ of solutions of equation \eqref{eq1}  for consecutive time levels $t_k$. 
The procedure of calculating the solution at the next time level $t_{k+1}$ is divided into two steps. First, on the interval $[t_k, t_{k+1}]$,  we solve a system of ODEs with the initial condition given by 
$\mu_{t_k} = \sum_{i=B}^{J}m^i_k\; \delta_{x^i_{k}},$
\begin{eqnarray}\label{b}
\mathrm{\frac{d}{dt}} x^i(t) = b_k(x^i(t)),
\quad x^i(t_k) = x_k^i\quad \mathrm{for}\;\; i= B, \dots, J,
\end{eqnarray}
where
$b_k(x) = b(t_{k}, \mu_{t_k})(x)$.
In the second step, we increase the index of each cohort by one, and then create a new empty cohort indexed by $B$ on the boundary. To proceed, we define a measure
$
\tilde \mu_{t_k}  = \sum_{i=B}^{J}m^i_k\; \delta_{x^i (t_{k+1})}$,
``freeze" coefficients $c$ and $\beta$ at the point $(t_k, \tilde \mu_{t_k})$, define
\begin{eqnarray}
\label{c_freeze}
c_k(x) = c \left( t_k, \tilde \mu_{t_k}\right)(x)\quad \mathrm{and}\quad
\beta_k(x) = \beta( t_k, \tilde \mu_{t_k})(x),
\end{eqnarray}
and  solve the following ODEs system on the time interval $t\in[t_k, t_{k+1}]$
\begin{eqnarray*}
\left\{
	\begin{array}{rcl}
\displaystyle
\mathrm{\frac{d}{dt}} m^i(t) &=& -c_k(x^i(t_{k+1})) m^i(t),\;\;\;
m^i(t_k) = m^i_k,
\quad\mathrm{for}\;\; i=B+1,\dots,J+1,
\\[2mm]
\displaystyle
\mathrm{\frac{d}{dt}} m^B(t) &=& -c_k(x^B(t_{k+1})) m^B(t) + \sum_{j=B}^{J+1} \beta_k(x^j(t_{k+1}))m^j(t),\;\;\;
m^B(t_k) = 0.
	\end{array}
\right.
\end{eqnarray*}
The output of a single step of the algorithm is thus 
$$
\mu_{t_{k+1}} = \sum_{i=B}^{J+1}m^i(t_{k+1}) \delta_{x^i(t_{k+1})}.
$$
The split-up algorithm is studied analytically and numerically in \cite{CGU}.

\section{Convergence of the algorithms}\label{Convergence}
\noindent
To estimate the distance between $\mu_t^n$ and the trajectory of semiflow $S$ starting at $\mu_0$ we use the following proposition, which allows us to consider equations locally in time.
\begin{proposition}\label{PropBressan}
Let $S: E \times [0, \delta] \times [0,T] \rightarrow E$ be a Lipschitz semiflow. For every Lipschitz continuous map $\nu:[0,T]\rightarrow E$ the following estimate holds,
\begin{equation}\label{PropBressan_formula}
\rho\big(\nu_t,S(t; 0)\mu_0\big) \leq L \int_{[0,t]} \liminf_{h \downarrow 0}\frac{\rho\big(\nu_{\tau+h}, S(h; \tau)\nu_{\tau}\big)}{h} d{\tau}, 
\end{equation}
where $\rho$ is a corresponding metric.
\end{proposition}
\noindent
The proof of Proposition \ref{PropBressan} is similar to the proof of Theorem~2.9 in \cite{bressan_00}. 
To apply Proposition \ref{PropBressan} we need to show the following result.
\begin{lemma}\label{lip_miu}
Let $\mu^n_t = \sum_{i} m^i(t) \delta_{x^i(t)}$, where $(x^i, m^i)$ is the output of the algorithm obtained by solving \eqref{EBT_org_system} with a boundary cohort defined as either in \eqref{EBT_org_BC} or \eqref
{EBT_simply_BC}. Then, 
$
\mu^n: [0, T] \to (\mathcal M^+ (\R_{\geq} ), \rho_F)
$
is Lipschitz continuous.
\end{lemma}
\begin{proofof}{Lemma \ref{lip_miu}}
Let $0 \leq s \leq t \leq T$ be such that $\modulo{t-s} \leq T^*$, where $T^*$ is the length of the interval of Lipschitz continuity of solutions of \eqref{EBT_org_system} with a boundary cohort defined as either in \eqref{EBT_org_BC} or \eqref
{EBT_simply_BC}. Without loss of generality we may assume that there is no internalization process on $(s,t)$.
By Lemma \ref{LemmaFlatMetric} we obtain
\begin{align*}
\rho_F(&\mu^n_t, \mu^n_s)
\leq
\sum_{i=B}^{J}\left(
m^i(s)\modulo{x^i(t) - x^i(s)} + \modulo{m^i(t) - m^i(s)}
\right)
\\
&\leq
({t-s})\sum_{i=B}^{J}\Big(
m^i(s)\mathrm{Lip}(x^i) + \mathrm{Lip}(m^i)
\Big)
\\
&\leq
({t-s})\max\{1,C\}\left(\sum_{i=B}^{J}
m^i(s) + \sum_{i=B}^{J}\mathrm{Lip}(m^i)\right).
\end{align*}
Due to the Theorem \ref{existence_org} and Theorem \ref{existence_org2} we know that $\sum_{i=B}^{J}m^i$ is bounded on each finite time interval. Additionally, using equations for $dm^i/dt$ and estimate \eqref{estimates_pi} from Lemma \ref{oszacowania} we conclude that $\sum_{i=B}^{J}\mathrm{Lip}(m^i)$ is also bounded independently of $J$ on each finite time interval.
\end{proofof}
\begin{theorem}\label{conv_1}
Let $\mu(\cdot)$ be a solution of \eqref{eq1} with the initial data $\mu_o$. Let $\mu_o^n = \sum_{i=B}^{N} m^i_o \delta_{x^i_o}$ be an approximation of $\mu_o$ by a sum of $n = N - B + 1$ Dirac deltas with the error of approximation
$$
e_x= \rho_F\left(\mu_o, \mu_o^n\right).
$$
Let
$\mu^n_t$ be the output of the EBT algorithm at time $t$ either for the original definition of the boundary cohort \eqref{EBT_org_BC} or the simplified one \eqref{EBT_simply_BC} with the initial condition $(x^i_o, m^i_o)$, $i = B, \dots, N$. Then, there exists a constant $C$ such that
$$
\rho_F(\mu(t), \mu^n_t) \leq e_x + C \Delta t.
$$
\end{theorem}
 \begin{remark}
 For a given Radon measure an arbitrarily good approximation in flat metric can be constructed in the form $\sum_{i=B}^{N} m^i_o \delta_{x^i_o}$ as shown in \cite{JMC}.
 \end{remark}
\begin{proofof}{Theorem \ref{conv_1}}
Let $\Delta t \leq T^*$, where $[0,T^*]$ is the time of existence of a unique solution of the EBT ODE system. 
Let $(\tau, \tau +h]$ contain no internalization time point.
As proved in \cite{CCGU}, problem \eqref{eq1} generates a Lipschitz semiflow $S$ such that $\mu(t) = S(t; 0)\mu_o$. In order to shorten the notation we define $\mu_{\tau+h} :=S(h; \tau) \mu^n_{\tau}$.
Measure  $\mu_{\tau +h}$ consists of $L =J-B+1$ Dirac deltas denoted henceforth as $n^i(\tau + t)\delta_{y^i(\tau + t)}$ and the density $f(t,\cdot)$ which arises due to the boundary condition.
The support of $f(t,\cdot)$ is contained in $[x_b, y^{abs}(t)]$, where $y^{abs}(\cdot)$ denotes the location of the characteristic line starting from $x_b$ at time $\tau$. 
We denote the total mass of $f(t,\cdot)$ as $n^{abs}(t) = \int_{x_b}^{y^{abs}(t)} f(t, x) d x$. For simplicity, in all
estimates below we will use a generic constant $C$ without specifying its exact form.

It is a classical result (see e.g. \cite{AmbrosioGigliSavare}) that a solution to the continuity equation is defined through so called ``push-forward" formula. Therefore,
\begin{equation}
\label{x_b_sim}
y^i(\tau + h) =  x^i(\tau) + \int_{\tau}^{\tau + h}b(t,\mu_t)(y^i(t))dt\quad \text{for}\;\; i= B, \dots, J.
\end{equation}
Equation on $y^{abs}$ is of the analogous form.
To obtain the expressions describing evolution of $n^i$ and $n^{abs}$ in time, we use proper test functions in the definition of weak solution \cite[Definition 2.2]{CCGU}. More precisely, let $\phi^i \in (C^1\cap W^{1,\infty})([0,T]\times \R_{\geq})$ be such that it is equal to $1$ on $(y^i(t) - \varepsilon, y^i(t)+\varepsilon)$ for $\varepsilon$ small enough and $\phi^i(t,y^j(t)) = 0$ for all $j \neq i$. Such a function exists, since the vector field $b$ is regular enough. This leads to
\begin{eqnarray}
\label{m^i_sim}
n^i(\tau + h) &=& m^i(\tau) - \int_{\tau}^{\tau + h}\!\!\!\!\!\!c( t, \mu_t)(y^i(t))n^i(t) dt.
\end{eqnarray}
Taking $\phi^{abs}  \in (C^1\cap W^{1,\infty})([0,T]\times \R_{\geq})$ such that $\phi^{abs}(t,x) = 1$ on $ [x_b, y^{abs(t)} + \varepsilon)$ and $\phi^i(t,y^j(t)) = 0$, $j = B, \dots, J$, leads to
\begin{eqnarray}
n^{abs}(\tau + h) &=& \sum_{i=B}^{J} \int_{\tau}^{\tau + h} 
\beta(t,\mu_t)(y^i(t)) n^i(t) dt + \mathcal{O}(h^2) .
\end{eqnarray}
We estimate
\begin{eqnarray}\label{estimate}
\rho_F(\mu_{\tau+h}, \mu^n_{\tau+h})&\leq& 
\sum_{i=B}^J m^i \left| x^i-y^i\right|+  \sum_{i=B+1}^J \left| n^i-m^i\right| \nonumber  \\&& +\left| n^B+n^{abs}-m^B\right| +n^{abs} \left| y^B-x_b\right|,
\end{eqnarray}
where all $x^i, y^i, m^i, n^i$ and $n^{abs}$ are evaluated at time $\tau + h$.
For the term $\modulo{x^i - y^i}$, $i=B+1, \dots, J$, we obtain
\begin{eqnarray}\label{estx}
&&
\left |x^i(\tau + h) - y^i(\tau + h)\right|
\leq
\int_{\tau}^{\tau + h} | b( t, \mu^n_t)(x^i(t)) -
b( t, \mu_t)(y^i(t)) |\; dt
\\
&& \leq 
\int_{\tau}^{\tau + h} \left(| b( t, \mu^n_t)(x^i(t)) -
b(t, \mu_t)(x^i(t)) | +
 | b(t, \mu_t)(x^i(t)) -
b(t, \mu_t)(y^i(t)) |\right) dt
\nonumber\\
&& \leq 
\norma{b}_{BC}  \left(\int_{\tau}^{\tau + h} \rho_F(\mu^n_t, \mu_t)\; dt +
 \int_{\tau}^{\tau + h} | x^i(t) - y^i(t)|\; dt \right)
\nonumber\\
&&\leq
\norma{b}_{BC}  \int_{\tau}^{\tau + h} \left(
\mathrm{Lip}(\mu^n)h + \rho_F(\mu^n_{\tau}, \mu_{\tau}) + \mathrm{Lip}(\mu)h
\right) dt 
\nonumber \\
&& \quad +\;
\norma{b}_{BC}  \int_{\tau}^{\tau + h} \left(
\mathrm{Lip}(x^i)h + \modulo{x^i(\tau) - y^i(\tau)}+ \mathrm{Lip}(y^i)h
\right) dt \leq C h^2,\nonumber
\end{eqnarray}
which holds due to the fact that $ \rho_F(\mu^n_{\tau}, \mu_{\tau}) = 0$ and $ \modulo{x^i(\tau) - y^i(\tau)} = 0$. 

We estimate the expression $m^B\modulo{x^B - y^B}$ by
\begin{eqnarray*}
\mathrm{Lip}(m^B) \Delta t \left ( \mathrm{Lip}(x^B) h + \mathrm{Lip}(y^B)h + \modulo{x^B(\tau) - y^B(\tau)}\right)
\leq Ch \Delta t.
\end{eqnarray*}
Furthermore, using boundedness of the total mass, we obtain  $$\sum_{i=B}^J m^i \left| x^i-y^i\right| \leq \sum_{i=B}^J m^i C h^2 =  \mathcal O(h^2) + \mathcal O(h)\Delta t.$$ 
Similar estimations as in \eqref{estx} using \eqref{m^i_sim} instead of \eqref{x_b_sim} lead to
\begin{eqnarray*}
&& \sum_{i=B}^J\modulo{m^i - n^i} \leq \norma{c}_{BC} \sum_{i=B}^J  n^i \int_{\tau}^{\tau + h} \left(
\mathrm{Lip}(\mu^n)h + \rho_F(\mu^n_{\tau}, \mu_{\tau}) + \mathrm{Lip}(\mu)h
\right) dt \\
&& \quad +\;
\norma{c}_{BC} \sum_{i=B}^J   \int_{\tau}^{\tau + h} \left(
\mathrm{Lip}(m^i)h + \modulo{m^i(\tau) - n^i(\tau)}+ \mathrm{Lip}(n^i)h
\right) dt = \mathcal O(h^2) \nonumber 
\end{eqnarray*}
due to boundedness of the total mass and boundedness of the sum of Lipschitz coefficients of $m^i$ and $n^i$.

Next, we estimate the term $\modulo{n^B+n^{abs}-m^{B}}$. For the simplified EBT algorithm, it holds
\begin{eqnarray*}
&&\modulo{n^B+n^{abs}-m^{B}} \leq \int_{\tau}^{\tau + h}
\modulo{
c(t, \mu^n_t)(x^B(t))m^B(t) - c( t, \mu_t)(y^B(t))n^B(t) } dt\\
&& +
\sum_{i=B}^{J}\int_{\tau}^{\tau + h}
\modulo{
\beta(t, \mu^n_t)(x^i(t))m^i(t) - \beta( t, \mu_t)(y^i(t))n^i(t) } dt + \mathcal O(h^2)=  \mathcal O(h^2).
\end{eqnarray*}
For the original definition of the boundary cohort we have a similar result, since the quantity $\modulo{\partial_x c(t,\mu)(x_b)\pi^B(t)} \leq C\Delta t$ for $t \in [0,\Delta t]$, which follows from Lemma \ref{oszacowania}. Finally, the last term
\begin{eqnarray*}
n^{abs}(\tau+h)\modulo{y^B (\tau+h)-x_b}\leq  C h \Delta t.
\end{eqnarray*}
Combining the above estimates, we obtain
\begin{eqnarray*}
\mathrm{liminf}_{h \to 0^+} \frac{1}{h}\rho_F(\mu_t^n, \mu_t^n)
& \leq& 
\mathrm{liminf}_{h \to 0^+}  \frac{1}{h} \left[
\mathcal{O}(h)\Delta t + \mathcal{O}(h^2) \right]
= C \Delta t,
\end{eqnarray*}
which, by Proposition \ref{PropBressan} implies that
$
\rho_F(\mu_t^n, \mu_t) \leq CL t \Delta t.
$
\end{proofof}
\begin{remark}
At the beginning of the proof we assumed that $\tau$ is not the internalization process. Otherwise, $n^B(t) = 0$ on $[\tau,\tau +h]$. However, the whole argumentation remains valid for this case.
\end{remark}
\noindent
Convergence of the split-up algorithm has been recently presented in \cite{CGU}.

\section{Numerical simulations}\label{Simulations}

The aim of this section is to confirm empirically the order of convergence and to compare accuracy and performance of the algorithms presented in this paper. 
\subsection{Test Case 1}
We consider a linear problem with the initial condition which is a stable stationary solution. The main goal of this test is to compare the influence of different approximations of the boundary cohort on the numerical solution. We choose the model parameters as in \cite{CGU}
$$
b(x) = 0.2(1 - x),\quad
c(x) = 0.2,\quad
\beta(x) = 2.4 (x^2 - x^3).
$$
The exact solution is $u(x,t) = 1$, $x \in [0,1]$.

\subsection{Test Case 2} We consider a nonlinear model, where the birth process depends on the total population state.
The aim of this test is to study influence of nonlocal terms on the numerical solution. Following \cite{kostova}, we set
\begin{eqnarray*}
b(x) &=& \mathrm{e}^{-x},
\quad
c(x)  = 
1 + \mathrm{e}^{-x} + \frac{\mathrm{e}^{-x}\sin(x)}{2 + \cos(x)},
\\
\beta(x, t, P) &=&
\frac{3}{2 + \cos(x)}
 \cdot
 \frac{0.5 + (1 + 0.5 \sin(1)) \mathrm{e}^{-t}}{0.5 + P},
\end{eqnarray*}
where $P = \int_{0}^{1}d \mu_t(x)$.
The exact solution is $u(x,t) = \mathrm{e}^{-t}(1 + 0.5\cos(x))$, $x \in (0,1]$.

\subsection{Test Case 3}\label{tc3} We consider an example with a large value of $\partial_{x}c(x_b)$. The main goal of this test is to show difficulties that may arise due to the lack of global in time existence of solutions to the original EBT scheme. The parameters are defined by
\begin{eqnarray*}
b(x) &=& \chi_{[0,0.5)} + (1-2(x-0.5))\chi_{[0.5,1]},
\\
c(x) &=&
\min(10, 10^4 x(1-x) ),
\\
\beta(x) &=&
10.
\end{eqnarray*}
The exact solution is not known analytically.

\subsection{Implementation of the algorithms}

The algorithms were implemented in C programming language in order to achieve maximum performance. All tests were performed on AMD Opteron(tm) Processor 8218 (each simulation running on a single core) and the binary executable was compiled using GCC 4.3.6 with -O3 optimization.
The implementation of the algorithms was based on solving systems of ODEs by the Euler explicit scheme. Derivatives of the functions $b$ and $c$, which are necessary for the original EBT algorithm, were computed symbolically. Therefore, cost and error of their numerical evaluation during the simulations were negligible. 

Errors of the schemes were calculated in the flat metric and $L_1$ norm. The flat metric was computed by the algorithm proposed in \cite{JMC}. To compute $L_1$ norm, we constructed a piecewise constant function, which bases on the numerical solution $\sum_{i=1}^{N}  m^i \delta_{x^i}$, as the following
\begin{equation}\label{density}
u(x,T) = \frac{2 m^i}{x_{i+1}-x_{i-1}}\;\;\; \mbox { for }\;\;\; x\in[\frac{x_{i-1}+x^i}{2},\frac{x^i+x_{i+1}}{2}).
\end{equation}
We compared such ``mollified" representation with the exact solutions in $L_1$ norm. The investigated particle methods do not guarantee convergence in such norms in general. Numerical results, however, show that for smooth parameters convergence can be expected.

In the third test case, where the exact solution is not known, we compared our results with the solution computed by the simplified EBT algorithm. The initial data was represented by 262144 nodes (Dirac deltas). During the simulation, we added 262144 boundary cohorts and did 16 steps of the Euler method between each internalization.

\subsection{Numerical results. Convergence of the algorithms and computational complexity}

In this section we compare three algorithms: Escalator Boxcar Train (EBT), Escalator Boxcar Train with simplified boundary conditions (sEBT) and Split-up (SU). We denote the number of initial nodes by $I$, the number of boundary cohorts created during a simulation by $K$ and the number of steps of the Euler method between each internalization by $J$. In most cases we set $I = KJ$.

Tables \ref{table1} - \ref{table4} show the numerical errors ($e_{I,K,J}$) and orders of convergence ($o_{I,K,J}$). The order of convergence is defined as  $o_{I,K,J} = \mathrm{log}_{2}( e_{I/2,K/2,J}/e_{I,K,J})$. All test cases are solved on the time interval $[0,1]$.

\begin{table*}[ht]
{\footnotesize
\begin{centering}
\begin{tabular}{|l|l l|l l|l l|}
\hline
&\multicolumn{2}{l|}{sEBT}&\multicolumn{2}{l|}{EBT}&\multicolumn{2}{l|}{SU} \\
I&Error&\textbf{Order}&Error&\textbf{Order}&Error&\textbf{Order} \\
\hline\hline
16 & 1.53e-02 & 1.03 & 1.31e-02 & 1.02 & 1.49e-02 & 1.04   \\ 
32 & 7.56e-03 & 1.02 & 6.56e-03 & 1.00 & 7.96e-03 & 0.90   \\ 
64 & 3.76e-03 & 1.01 & 3.28e-03 & 1.00 & 4.14e-03 & 0.94   \\ 
128 & 1.88e-03 & 1.00 & 1.64e-03 & 1.00 & 2.11e-03 & 0.97   \\ 
256 & 9.36e-04 & 1.00 & 8.20e-04 & 1.00 & 1.07e-03 & 0.99   \\ 
512 & 4.68e-04 & 1.00 & 4.10e-04 & 1.00 & 5.36e-04 & 0.99   \\ 
1024 & 2.34e-04 & 1.00 & 2.05e-04 & 1.00 & 2.68e-04 & 1.00   \\ 
2048 & 1.17e-04 & 1.00 & 1.03e-04 & 1.00 & 1.34e-04 & 1.00   \\ 
4096 & 5.84e-05 & 1.00 & 5.13e-05 & 1.00 & 6.73e-05 & 1.00   \\ 
8192 & 2.92e-05 & 1.00 & 2.56e-05 & 1.00 & 3.36e-05 & 1.00   \\ 
16384 & 1.46e-05 & 1.00 & 1.28e-05 & 1.00 & 1.68e-05 & 1.00   \\ 
32768 & 7.30e-06 & 1.00 & 6.41e-06 & 1.00 & 8.41e-06 & 1.00   \\ 
65536 & 3.65e-06 & 1.00 & 3.20e-06 & 1.00 & 4.21e-06 & 1.00   \\ 
131072 & 1.83e-06 & 1.00 & 1.60e-06 & 1.00 & 2.10e-06 & 1.00   \\ 
262144 & 9.13e-07 & 1.00 & 8.01e-07 & 1.00 & 1.05e-06 & 1.00   \\ 
524288 & 4.56e-07 & 1.00 & 4.01e-07 & 1.00 & 5.26e-07 & 1.00   \\ 
1048576 & 2.28e-07 & 1.00 & 2.00e-07 & 1.00 & 2.63e-07 & 1.00   \\ 
\hline
\end{tabular}
\par\end{centering}
\vspace{0.3cm}
\caption{\footnotesize Test Case 1. Numerical error and order of convergence measured in flat metric, $J=4$, $K=I/4$.}\label{table1}
}
\end{table*}

\begin{table}[ht]
{\footnotesize
\begin{centering}
\begin{tabular}{|l|l l|l l|l l|}
\hline
 &\multicolumn{2}{l|}{sEBT}&\multicolumn{2}{l|}{EBT}&\multicolumn{2}{l|}{SU} \\
I&Error&\textbf{Order}&Error&\textbf{Order}&Error&\textbf{Order} \\
\hline\hline
16 & 4.68e-02 & 0.88 & 1.09e-02 & 0.73 & 9.99e-02 & 0.98   \\ 
32 & 2.43e-02 & 0.94 & 5.89e-03 & 0.89 & 4.98e-02 & 1.01   \\ 
64 & 1.24e-02 & 0.97 & 3.05e-03 & 0.95 & 2.48e-02 & 1.01   \\ 
128 & 6.26e-03 & 0.99 & 1.55e-03 & 0.98 & 1.23e-02 & 1.00   \\ 
256 & 3.14e-03 & 0.99 & 7.83e-04 & 0.99 & 6.15e-03 & 1.00   \\ 
512 & 1.57e-03 & 1.00 & 3.93e-04 & 0.99 & 3.07e-03 & 1.00   \\ 
1024 & 7.88e-04 & 1.00 & 1.97e-04 & 1.00 & 1.54e-03 & 1.00   \\ 
2048 & 3.94e-04 & 1.00 & 9.86e-05 & 1.00 & 7.68e-04 & 1.00   \\ 
4096 & 1.97e-04 & 1.00 & 4.93e-05 & 1.00 & 3.84e-04 & 1.00   \\ 
8192 & 9.86e-05 & 1.00 & 2.47e-05 & 1.00 & 1.92e-04 & 1.00   \\ 
16384 & 4.93e-05 & 1.00 & 1.23e-05 & 1.00 & 9.60e-05 & 1.00   \\ 
32768 & 2.47e-05 & 1.00 & 6.17e-06 & 1.00 & 4.80e-05 & 1.00   \\ 
65536 & 1.23e-05 & 1.00 & 3.08e-06 & 1.00 & 2.40e-05 & 1.00   \\ 
131072 & 6.17e-06 & 1.00 & 1.54e-06 & 1.00 & 1.20e-05 & 1.00   \\ 
262144 & 3.08e-06 & 1.00 & 7.72e-07 & 1.00 & 6.00e-06 & 1.00   \\ 
524288 & 1.55e-06 & 1.00 & 3.90e-07 & 0.99 & 3.00e-06 & 1.00   \\ 
1048576 & 7.84e-07 & 0.98 & 2.06e-07 & 0.92 & 1.51e-06 & 0.99   \\ 
\hline
\end{tabular}
\par\end{centering}
\vspace{0.3cm}
\caption{\footnotesize Test Case 1. Numerical error and order of convergence measured by L1 norm, $J=4$, $K=I/4$.}
}
\end{table}

\begin{table}[ht]
{\footnotesize
\begin{centering}
\begin{tabular}{|l|l l|l l|l l|}
\hline
 &\multicolumn{2}{l|}{sEBT}&\multicolumn{2}{l|}{EBT}&\multicolumn{2}{l|}{SU} \\
I&Error&\textbf{Order}&Error&\textbf{Order}&Error&\textbf{Order} \\
\hline\hline
16 & 6.09e-02 & 1.12 & 6.25e-02 & 1.12 & 1.29e-01 & 0.82   \\ 
32 & 3.67e-02 & 0.73 & 3.92e-02 & 0.67 & 5.72e-02 & 1.17   \\ 
64 & 1.63e-02 & 1.17 & 1.72e-02 & 1.19 & 3.06e-02 & 0.90   \\ 
128 & 9.32e-03 & 0.81 & 1.01e-02 & 0.77 & 1.40e-02 & 1.13   \\ 
256 & 5.02e-03 & 0.89 & 5.41e-03 & 0.90 & 6.78e-03 & 1.04   \\ 
512 & 2.27e-03 & 1.15 & 2.46e-03 & 1.14 & 3.52e-03 & 0.95   \\ 
1024 & 1.19e-03 & 0.93 & 1.29e-03 & 0.94 & 1.72e-03 & 1.03   \\ 
2048 & 6.37e-04 & 0.90 & 6.87e-04 & 0.91 & 8.42e-04 & 1.03   \\ 
4096 & 2.92e-04 & 1.12 & 3.18e-04 & 1.11 & 4.33e-04 & 0.96   \\ 
8192 & 1.56e-04 & 0.91 & 1.69e-04 & 0.91 & 2.12e-04 & 1.03   \\ 
16384 & 6.97e-05 & 1.16 & 7.59e-05 & 1.15 & 1.11e-04 & 0.94   \\ 
32768 & 3.54e-05 & 0.98 & 3.85e-05 & 0.98 & 5.48e-05 & 1.01   \\ 
65536 & 1.83e-05 & 0.95 & 1.99e-05 & 0.96 & 2.70e-05 & 1.02   \\ 
131072 & 9.74e-06 & 0.91 & 1.05e-05 & 0.91 & 1.32e-05 & 1.03   \\ 
262144 & 4.35e-06 & 1.16 & 4.74e-06 & 1.15 & 6.91e-06 & 0.94   \\ 
\hline
\end{tabular}
\par\end{centering}
\vspace{0.3cm}
\caption{\footnotesize Test Case 2. Numerical error and order of convergence measured by flat metric, $J=4$, $K=I/4$.}
}
\end{table}

\begin{table*}
{\footnotesize
\begin{centering}
\begin{tabular}{|l|l l|l l|l l|}
\hline
 &\multicolumn{2}{l|}{sEBT}&\multicolumn{2}{l|}{EBT}&\multicolumn{2}{l|}{SU} \\
I&Error&\textbf{Order}&Error&\textbf{Order}&Error&\textbf{Order} \\
\hline\hline
16 & 9.90e-02 & 0.87 & 6.39e-02 & 0.82 & 1.68e-01 & 0.99   \\ 
32 & 5.71e-02 & 0.79 & 4.13e-02 & 0.63 & 8.34e-02 & 1.01   \\ 
64 & 2.86e-02 & 1.00 & 2.16e-02 & 0.94 & 4.08e-02 & 1.03   \\ 
128 & 1.52e-02 & 0.91 & 1.15e-02 & 0.91 & 2.08e-02 & 0.97   \\ 
256 & 7.98e-03 & 0.93 & 6.21e-03 & 0.89 & 1.07e-02 & 0.95   \\ 
512 & 3.80e-03 & 1.07 & 2.88e-03 & 1.11 & 5.13e-03 & 1.06   \\ 
1024 & 1.95e-03 & 0.97 & 1.50e-03 & 0.94 & 2.60e-03 & 0.98   \\ 
2048 & 1.02e-03 & 0.93 & 7.95e-04 & 0.92 & 1.35e-03 & 0.95   \\ 
4096 & 4.85e-04 & 1.08 & 3.72e-04 & 1.10 & 6.48e-04 & 1.06   \\ 
8192 & 2.52e-04 & 0.94 & 1.96e-04 & 0.93 & 3.34e-04 & 0.96   \\ 
16384 & 1.19e-04 & 1.08 & 9.07e-05 & 1.11 & 1.60e-04 & 1.06   \\ 
32768 & 5.95e-05 & 1.00 & 4.53e-05 & 1.00 & 7.98e-05 & 1.00   \\ 
65536 & 3.03e-05 & 0.97 & 2.33e-05 & 0.96 & 4.05e-05 & 0.98   \\ 
131072 & 1.58e-05 & 0.94 & 1.22e-05 & 0.93 & 2.09e-05 & 0.96   \\ 
262144 & 7.44e-06 & 1.08 & 5.67e-06 & 1.11 & 9.98e-06 & 1.06   \\ 
\hline
\end{tabular}
\par\end{centering}
\vspace{0.3cm}
\caption{\footnotesize Test Case 2. Numerical error and order of convergence measured by L1 norm, $J=4$, $K=I/4$.}\label{table4}
}
\end{table*}

\begin{table}
{\footnotesize
\begin{centering}
\begin{tabular}{|l l|l|l|l|l|}
\hline
   & &sEBT&EBT&SU&EBT* \\
I=K&J&Error&Error&Error&Error \\
\hline\hline
32 & 1 & 6.41e-02 & 13.36 & 7.75e-01 & 13.36\\ 
32 & 2 & 6.39e-02 & 1.11 &  2.08e-01 & 1.11 \\ 
32 & 8 & 6.51-02 & $x^B > x^1$ & 6.25e-01 & 1.05 \\ 
32 & 32 & 6.59e-02 & $m^B<0$ & 8.98e-01 & 1.05\\ 
128 & 1 &   5.59e-02 & 1.10 & 1.31e-01 & 1.10 \\ 
128 & 2 &   5.59e-02 & 1.63e-01 & 9.09e-02 & 1.63-01 \\ 
128 & 8 &   5.56e-02 & 3.34e-01 & 2.77e-01 & 3.24-01 \\ 
128 & 32 &  5.39e-02 & 4.72e-01 & 3.27e-01 & 4.52-01 \\ 
1024 & 1 &  4.97e-02 & 5.36e-02 & 5.36e-02 & 5.36-02 \\ 
1024 & 2 &  3.83e-02 & 2.57e-02 & 2.87e-02 & 2.57-02 \\ 
1024 & 8 &  3.49e-02 & 2.21e-03 & 9.76e-03 & 2.24-03 \\ 
1024 & 32 & 3.48e-02 & 3.41e-03 & 5.00e-03 & 3.37-03 \\ 
\hline
\end{tabular}
\par\end{centering}
\vspace{0.3cm}
\caption{\footnotesize Test Case 3. Numerical error measured in flat metric. In this test case we assumed that the number of initial nodes is equal to the number of boundary cohorts. We compared the influence of solving the corresponding ODEs systems with increasing precision, i.e. imposing more Euler steps between each internalization. This example also shows problems which may arise due to the lack of global in time existence of solutions to the original EBT ODEs system. }\label{Table5}
}
\end{table}

\begin{figure}
\begin{centering}
\includegraphics[scale=0.9]{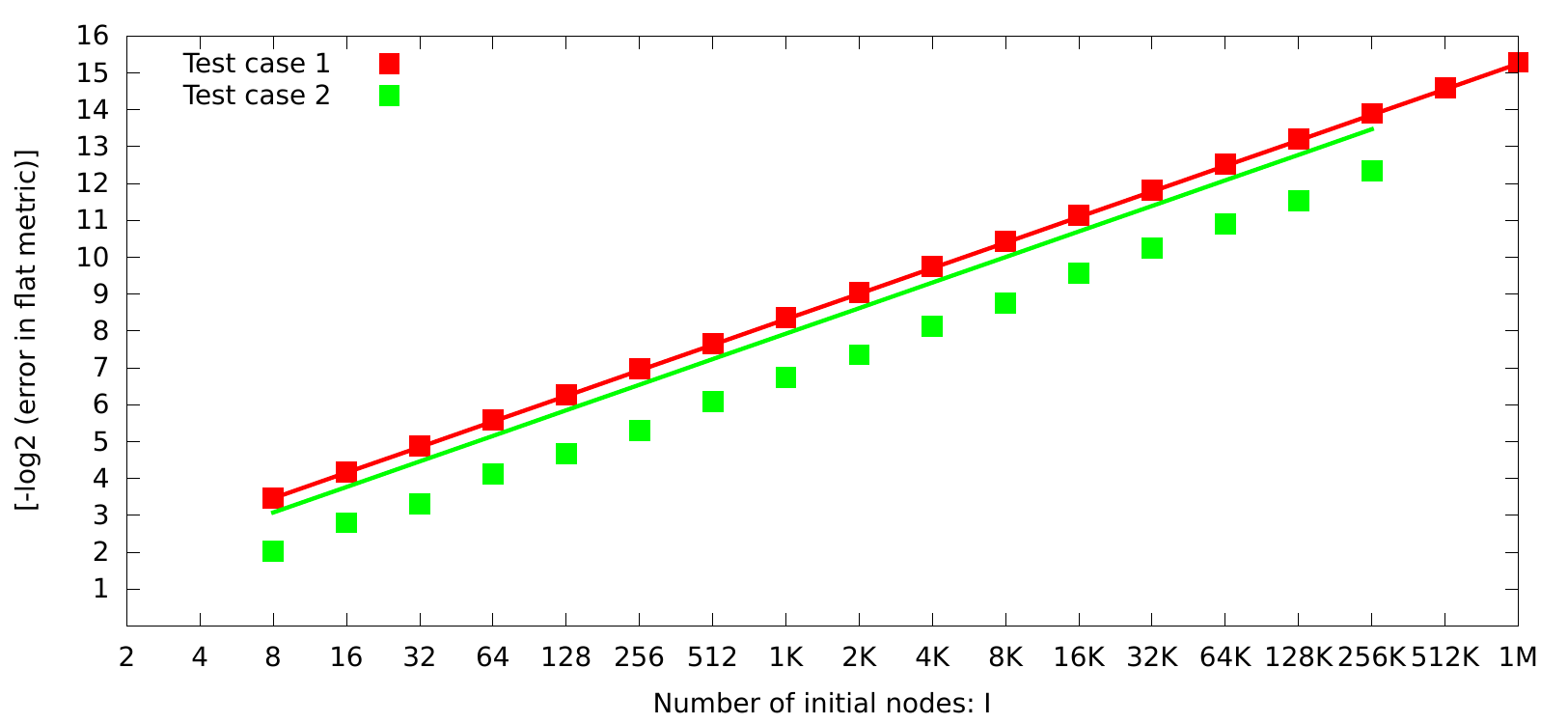}
\caption{\footnotesize The error of sEBT scheme in flat metric for $J=4,K=I/4$. 
Red dots represent Test Case 1, green dots - Test Case 2. Solid lines represent errors of approximation of the initial conditions in the respective test cases.}
\end{centering}
\end{figure}

\begin{figure}
\begin{centering}
\includegraphics[scale=0.9]{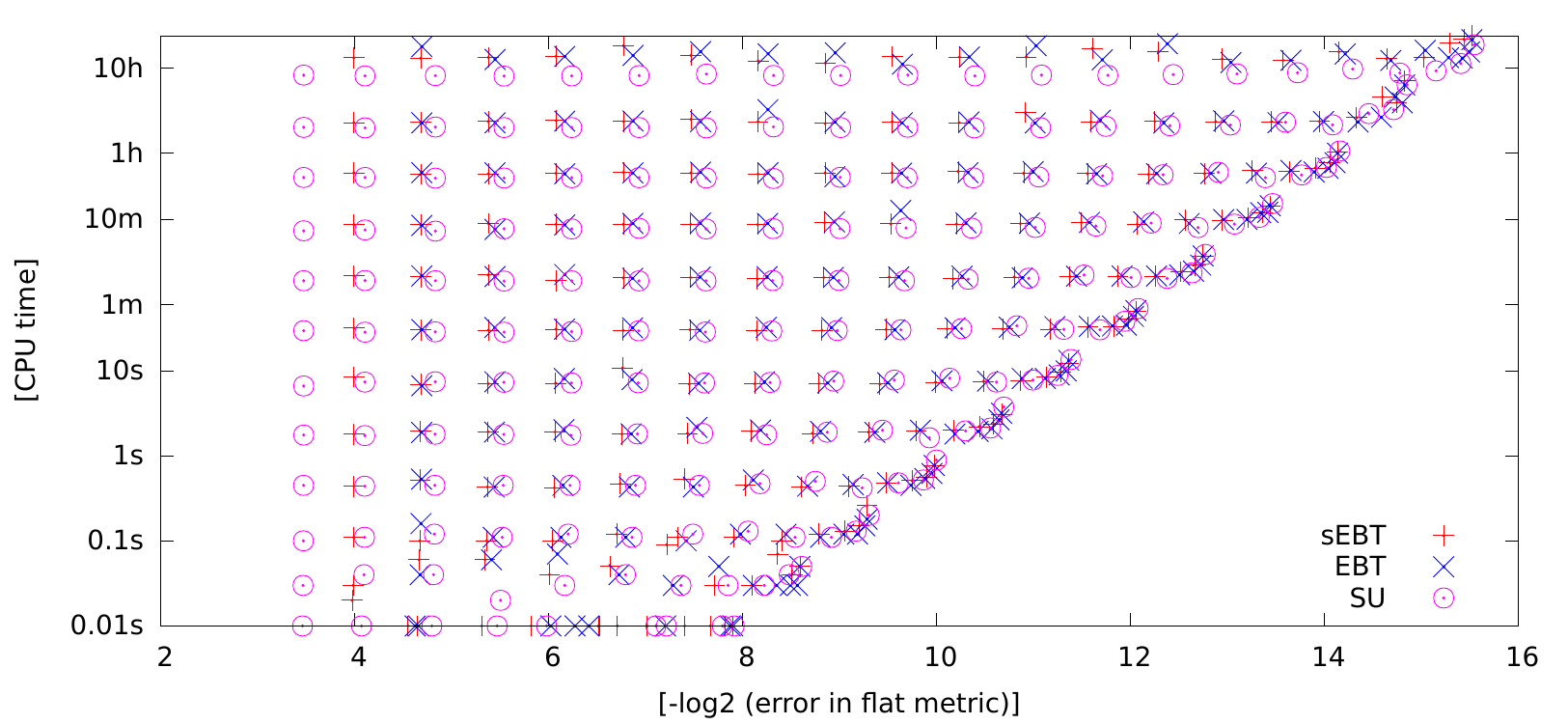}
\caption{\footnotesize CPU time required to achieve given accuracy. 
Each point in this plot represents a simulation for $I=2^i$, $K=2^k$, $K J=I$ and $i,k\in\{2,3,..,19\}$.
The points furthest to the right (high accuracy) correspond to these solutions for which $K$ was close to $I$. No significant difference in efficiency between algorithms could be found. The plot also suggests that maximum efficiency for this problem can be achieved by choosing $J\leq 4$.}
\end{centering}
\end{figure}

\begin{figure}
\begin{centering}
\includegraphics[scale=0.8]{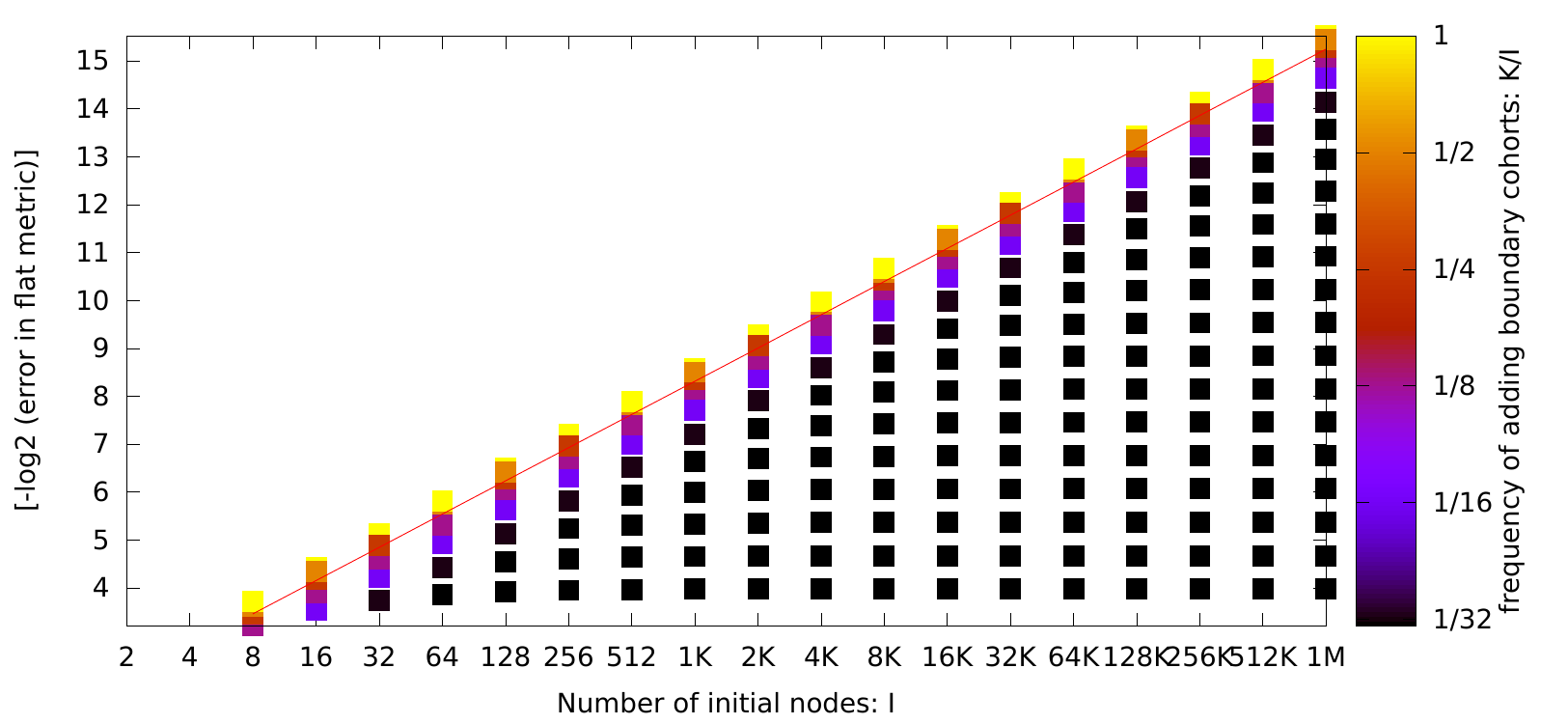}
\includegraphics[scale=0.8]{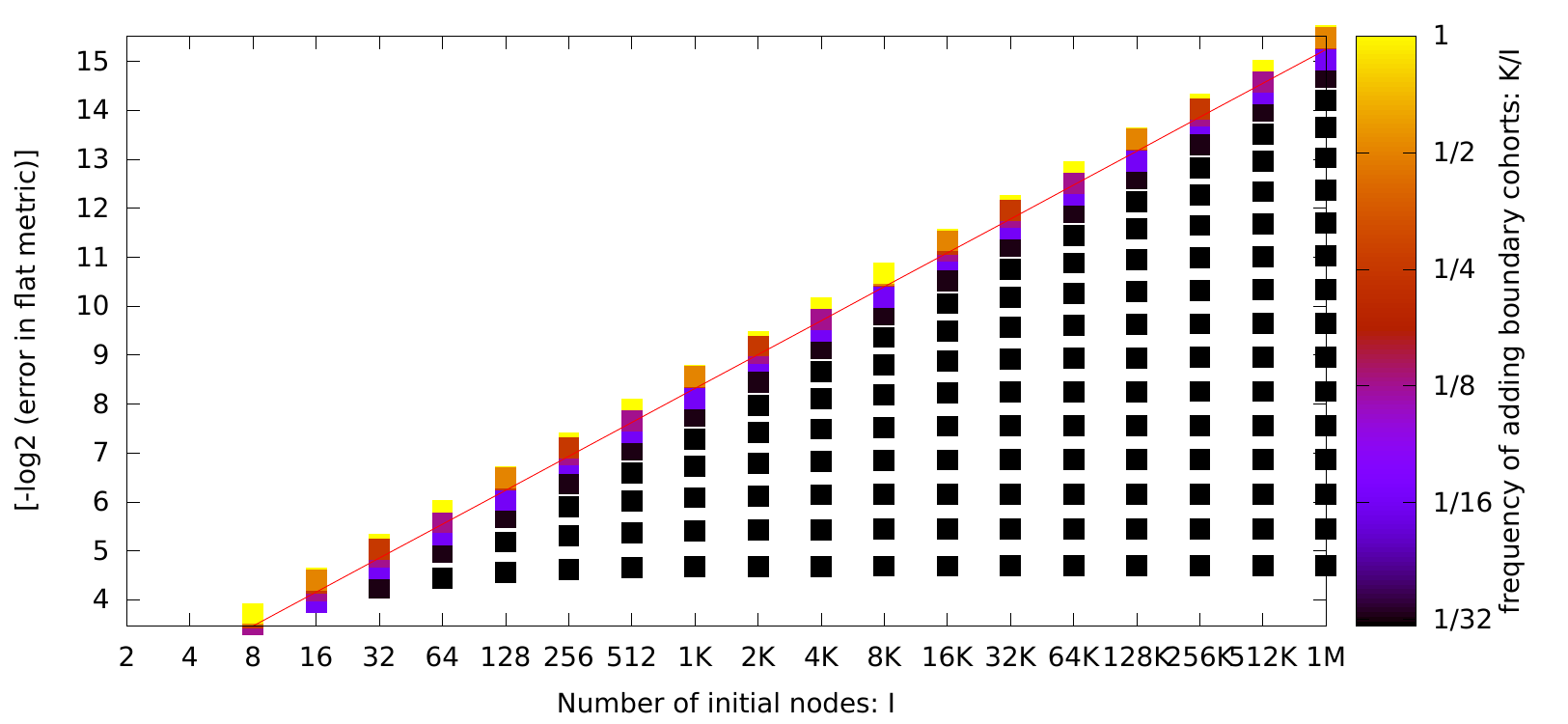}
\includegraphics[scale=0.8]{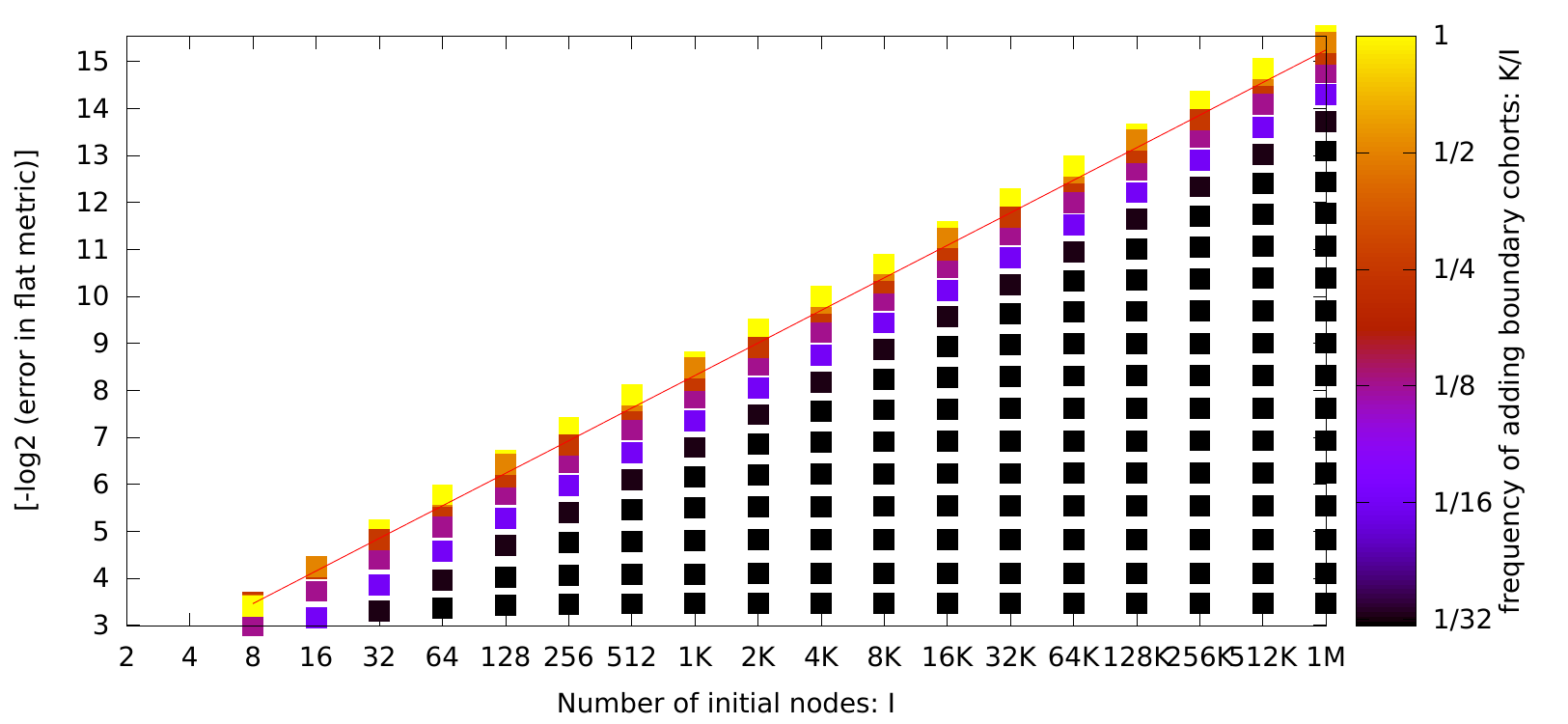}
\caption{\footnotesize Full map of errors for Test Case 1 (from above: sEBT, EBT, SU). The plots show the dependence of numerical error in flat metric upon number of initial nodes $I$ (X axis) and frequency of adding boundary cohorts $K/I = 1/J$ (color) for Test Case 1. The solid line represents the accuracy of initial condition approximation.}
\end{centering}
\end{figure}

\begin{figure}
\begin{centering}
\includegraphics[scale=0.8]{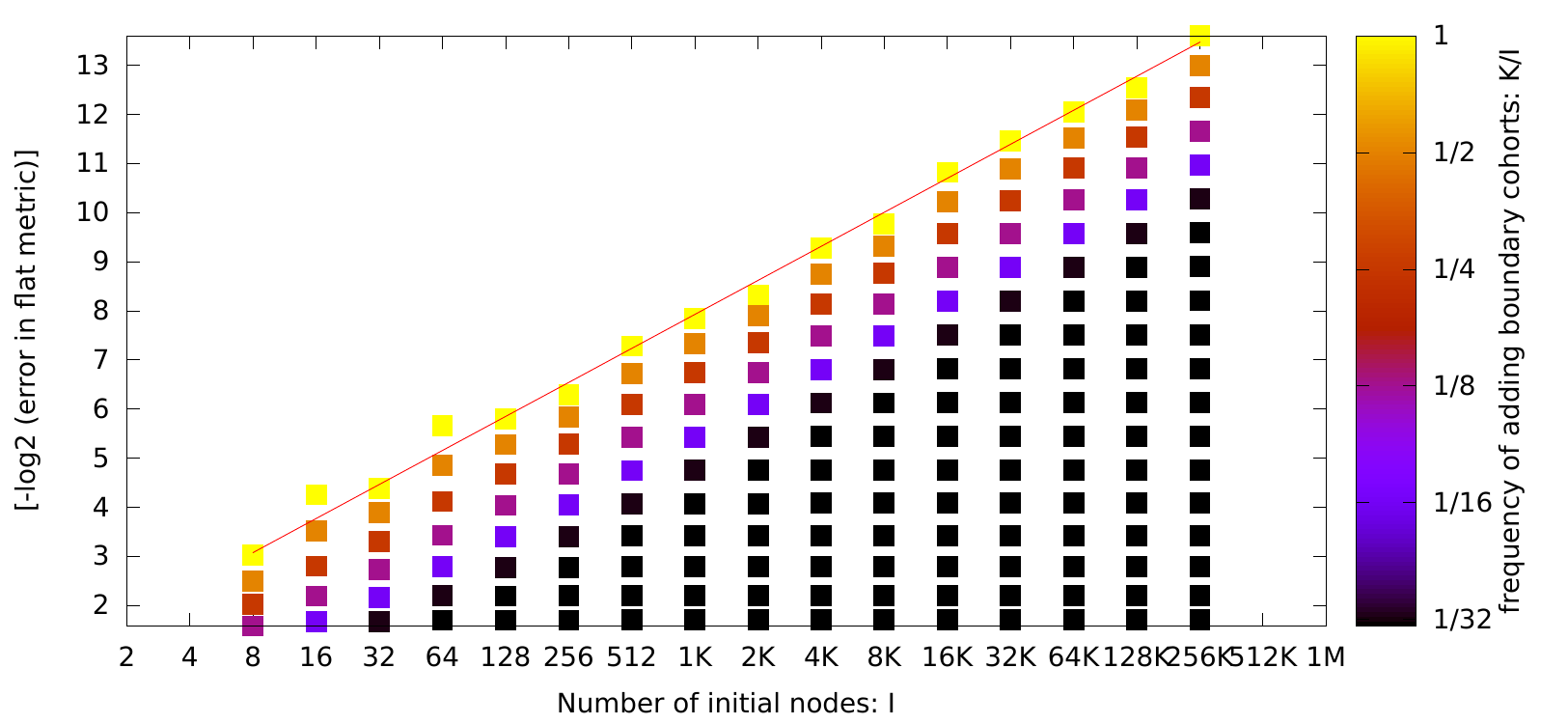}
\includegraphics[scale=0.8]{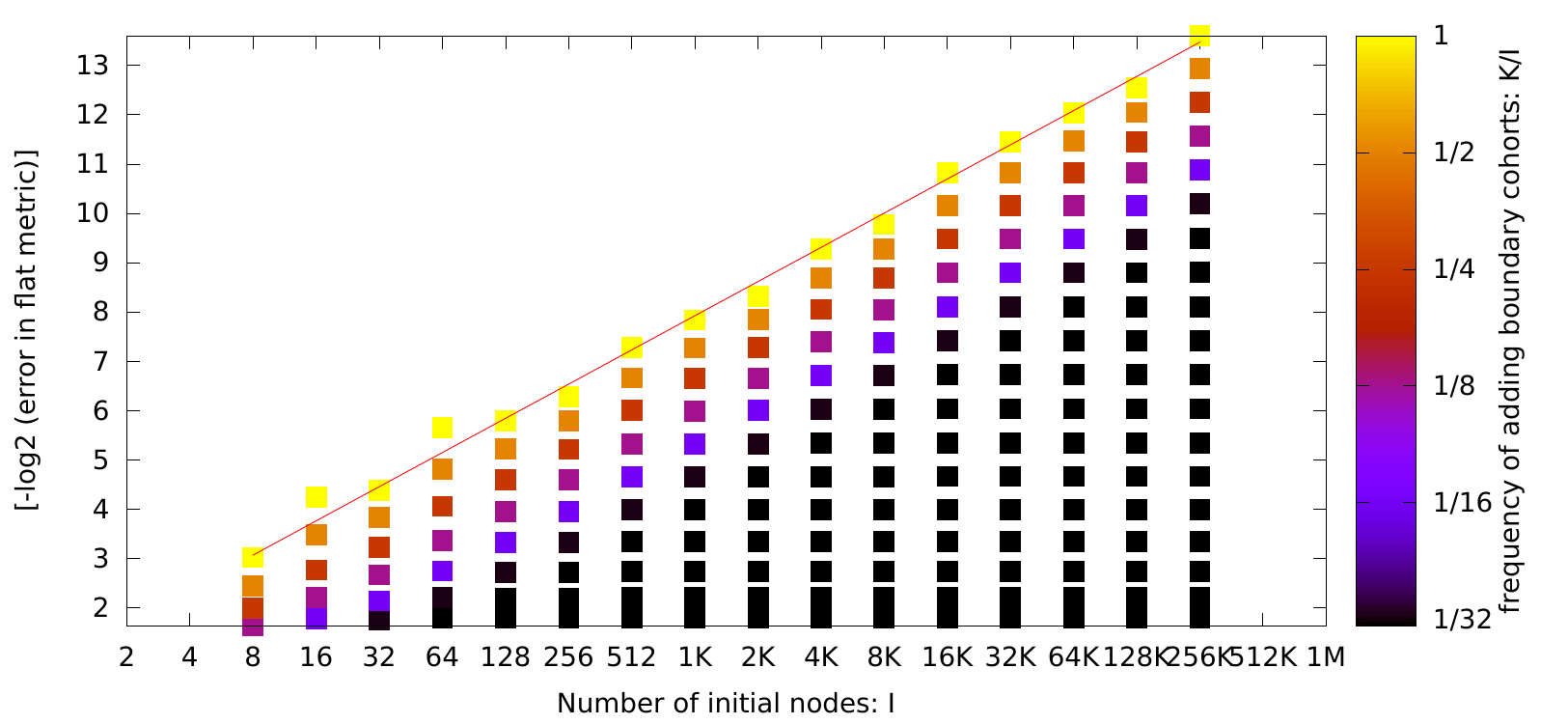}
\includegraphics[scale=0.8]{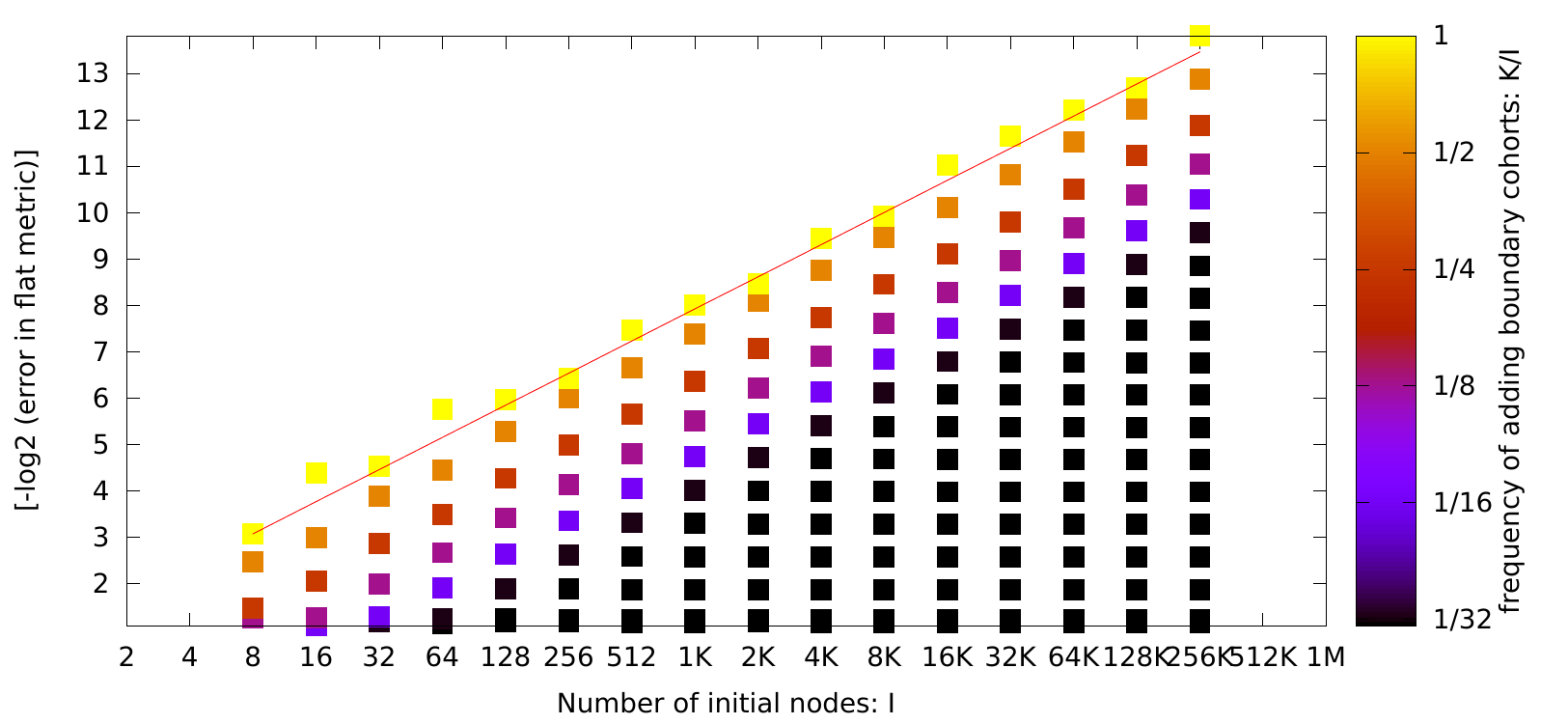}
\caption{\footnotesize Full map of errors for Test Case 2 (from above: sEBT, EBT, SU). The plots show the dependence of numerical error in flat metric upon number of initial nodes $I$ (X axis) and frequency of adding boundary cohorts $K/I = 1/J$ (color) for Test Case 2. The solid line represents the accuracy of initial condition approximation.}
\end{centering}
\end{figure}

\subsection{Conclusions}
Numerical simulations confirmed linear convergence of the presented algorithms. The optimal choice of the parameters $I$ and $J$, however, greatly depends on the specific model, which has to be solved. All three algorithms share the same complexity $O((I+K)KJ)$, therefore there are no major differences in performance. 

In majority of tests, the SU is $15\%-20\%$ faster than the other algorithms, most likely due to the sequential operations on consecutive elements of tables. Usually, it is also the least accurate. Accuracy of the sEBT and EBT depends greatly on the behaviour of functions $b$ and $c$ around $x_b = 0$. In the case of steep growth or decline, the EBT gives better results. The performance of these two does not differ by more than $5\%$ in all test cases.

For the low values of $K$, the original EBT may reveal some unexpected behavior (as shown in Table \ref{Table5}), as there is no guarantee that $x^B$ remains smaller than $x^1$ globally in time, and even that $m^B$ remains non-negative. In models where $\partial_{x}c(x_b)$ and $\partial_{x}b(x_b)$ are large this can be a serious problems, as $K$ may need to be extremly large to prevent this unwanted behavior.

\section{Appendix}

\subsection{Derivation of the EBT scheme}
To streamline the presented analysis, we focus on the linear case and assume that a solution $u(t,\cdot)$ is a compactly supported and integrable function, which leads to the following problem
\begin{eqnarray}\label{ebt_toy}
\partial_t u(t,x) + \partial_x (b(t,x)u(t,x)) + c(t,x)u(t,x) &=& 0
\\
\nonumber
b(t, x_b)u(t,x_b) &=& \int_{x_b}^{+\infty}\beta(t,y)u(t,u) d y. \end{eqnarray}
We also require that the model functions $b$, $c$ and $\beta$ are more regular, i.e.
$$
b,c,\beta : [0,T] \to C^2([x_b, +\infty)).
$$ 
Higher regularity is imposed in order to apply the first order Taylor approximation.
Let $\{\Omega_i(t)\}_{i=B}^{J}$ be a collection of pairwise disjoint intervals
$
\Omega_i(t) = [l_i(t), l_{i+1}(t)),
$
such that
$
\mathrm{supp}(u_o) \subset \bigcup_{i=0}^M \Omega_i(t) =: \Omega(t),
$
where $u_o$ is the initial distribution of individuals.
Set $\Omega_i(t)$ denotes a range of $i$-th cohort. 
Boundaries between the cohorts evolve according to the equation
\begin{equation*}
\mathrm{\frac{d}{dt}}l_i(t) = b(t, u)(l_i(t)), \quad
\mathrm{for}\;\;\; i=B, \dots, J.
\end{equation*}
A lower bound of the boundary cohort is constant in time,  $l_B(t) = x_b$.  The number of individuals is given by
$$
m^i(t) = \int_{\Omega_i(t)} u(t,x)dx
$$ and the average value of the structure variable within the cohort by
\begin{equation}\label{x^i}
x^i(t) = \frac{1}{m^i(t)} \int_{\Omega_i(t)} x u(t,x)\; dx,\quad
i=B, \dots, J.
\end{equation}
Since the boundary cohort is initially empty, we define 
$$
\pi^B(t) = \int_{\Omega_B(t)}(x - x_b)u(t,x) dx
\quad\mathrm{and}\quad
x^B(t) = x_b + \frac{\pi^B(t)}{m^B(t)}.
$$
Differentiating $m^i$, $x^i$ and $\pi^B$ yields  for $i=B+1, \dots, J$ 
\begin{eqnarray*}
\mathrm{\frac{d}{dt}} m^i(t) &=&
 \int_{\Omega_i(t)} \partial_t u(t,x) dx + \mathrm{\frac{d}{dt}} l_{i+1}(t) u(t, l_{i+1}(t)) - \mathrm{\frac{d}{dt}} l_i(t) u(t, l_i(t))
\\
 &=&
 \int_{\Omega_i(t)} \partial_t u(t,x) dx + b(t, l_{i+1}(t)) u(t, l_{i+1}(t)) - b(t, l_i(t)) u(t, l_i(t))
\\
 &=&
 \int_{\Omega_i(t)} \partial_t u(t,x) dx +
\int_{\Omega_i(t)} \partial_x \left(
b(t,x)u(t,x)
\right)
=
- \int_{\Omega_i(t)} c(t,x) u(t,x) dx.
\end{eqnarray*}
Similarly, for the boundary cohort it holds that
\begin{eqnarray*}
\mathrm{\frac{d}{dt}} m^B(t) &=&
\int_{\Omega_B(t)}\partial_t u(t,x) dx + b(t,l_1(t))u(t, l_1(t)) -
b(t,x_b)u(t,x_b) + b(t,x_b)u(t,x_b)
\\
&=&
-\int_{\Omega_B(t)}c(t,x) u(t,x)  dx+ \int_{x_b}^{+\infty}\beta(t,u)u(t,y) \d y.
\end{eqnarray*}
We calculate first moments of the internal cohorts $y(t) = \int_{\Omega_i(t)} x u(t,x) dx$,
\begin{eqnarray*}
\mathrm{\frac{d}{dt}} y_i(t) &=&
\int_{\Omega_i(t)}x \; \partial_tu(t,x) dx +
b(t,l_{i+1}(t))u(t, l_{i+1}(t)) \;  l_{i+1}(t) -
b(t,l_{i}(t))u(t, l_{i}(t)) \;  l_{i}(t) 
\\
&=&
\int_{\Omega_i(t)}x \; \partial_tu(t,x) dx +
\int_{\Omega_i(t)} \partial_x\Big(x (b(t,x)u(t,x) \Big) dx 
\\
&=&
\int_{\Omega_i(t)}x \; \partial_tu(t,x) dx +
\int_{\Omega_i(t)} x \partial_x (b(t,x)u(t,x)) dx +
\int_{\Omega_i(t)} b(t,x)u(t,x) dx 
\\
&=&
-\int_{\Omega_i(t)}x c(t,x)u(t,x) dx +
\int_{\Omega_i(t)} b(t,x)u(t,x) dx.
\end{eqnarray*}

\begin{eqnarray*}
\mathrm{\frac{d}{dt}}x^i(t) &=& \frac{\mathrm{\frac{d}{dt}}y^i(t)}{m^i(t)} -
\frac{y^i(t) \mathrm{\frac{d}{dt}}m^i(t)}{m^i(t)^2}
=
 \frac{\mathrm{\frac{d}{dt}}y^i(t)}{m^i(t)} -
x^i(t)\mathrm{\frac{d}{dt}}m^i(t)
\\
&=&
- \;\frac{1}{m^i(t)} \int_{\Omega_i(t)} x c(t, x) u(t,x) dx
+
\frac{1}{m^i(t)} \int_{\Omega_i(t)}
 b(t, x)u(t,x) dx
\\
&&\quad +\;
\frac{1}{m^i(t)} \int_{\Omega_i(t)} x^i(t) c(t,x)u(t,x)dx
\\
&=&
\frac{1}{m^i(t)} \int_{\Omega_i(t)}(x^i(t) - x) c(t, x) u(t,x) dx
+
\frac{1}{m^i(t)} \int_{\Omega_i(t)}
 b(t, x)u(t,x) dx.
\end{eqnarray*}
Finally, for the boundary cohort we obtain
\begin{eqnarray*}
\mathrm{\frac{d}{dt}}\pi^B(t) &=& 
\mathrm{\frac{d}{dt}}\left(
\int_{\Omega_B(t)}x u(t,x)dx
\right)
-
x_b \mathrm{\frac{d}{dt}}m^B(t)
\\
&=&
\int_{\Omega_B(t)}x \; \partial_tu(t,x) dx +
b(t,l_{1}(t))u(t, l_{1}(t)) \;  l_{1}(t) \pm
b(t,x_b)u(t, x_b) \;  x_b -
x_b \mathrm{\frac{d}{dt}}m^B(t)
\\
&=&
-\int_{\Omega_B(t)}x c(t,x)u(t,x) dx +
\int_{\Omega_B(t)} b(t,x)u(t,x) dx  + x_b\int_{x_b}^{+\infty}\beta(t,y)u(t,y) \d y 
\\
&&\;
+ x_b\int_{\Omega_B(t)}c(t,x) u(t,x)  -
x_b \int_{x_b}^{+\infty}\beta(t,u)u(t,y) \d y
\\
&=&
\int_{\Omega_B(t)}(x_b-x) c(t,x)u(t,x) dx +
\int_{\Omega_B(t)} b(t,x)u(t,x) dx .
\end{eqnarray*}
\noindent
\textbf{Approximation.}
To obtain a closed form of the scheme, we approximate $b$, $c$ and $\beta$. Using \eqref{x^i}, we obtain
\begin{equation}\label{mean}
\begin{split}
 \int_{\Omega_i(t)} (x^i(t) - x) u(t,x)dx = x^i(t) m^i(t) - x^i(t)m^i(t) = 0.
\end{split}
\end{equation}
Moreover,  it holds for $f \in C^2$
\begin{eqnarray*}
\int_{\Omega(t)} f(x)  u(t,x)dx
 &=& 
\sum_{i=B}^{J} \int_{\Omega_i(t)} f(x^i(t)) u(t,x)dx+
\sum_{i=B}^{J}\int_{\Omega_i(t)} \mathrm{\frac{d}{dx}}f(x^i(t))(x - x^i(t)) u(t,x)dx
\\
&& + \sum_{i=B}^{J} \int_{\Omega_i(t)} \mathcal O(\modulo{x - x^i(t)}^2) u(t,x)dx
=
\sum_{i=B}^{J}f(x^i(t)) m^i(t) +
\mathcal O( t^2).
\end{eqnarray*}
Therefore, the first order approximation is given as
\begin{eqnarray}\label{Taylor_app}
\int_{\Omega(t)} f(x)  u(t,x)dx
&=&
\sum_{i=B}^{J}f(x^i(t)) m^i(t).
\end{eqnarray}
\noindent
Application of \eqref{Taylor_app} and neglecting the second (and higher) order terms  yields
\begin{eqnarray*}
\mathrm{\frac{d}{dt}}m^i(t)
&=& c(t, x^i(t)) m^i(t),
\\
\mathrm{\frac{d}{dt}}x^i(t) &=&
 b(t, x^i(t)).
\end{eqnarray*}
For the boundary cohort we expand $c(x,t)$ around $x_b$
\begin{eqnarray*}
\mathrm{\frac{d}{dt}}m^B(t) &\!\!\!=& \!\!\!
-\int_{\Omega_B(t)}c(t,x) u(t,x) dx + \int_{x_b}^{+\infty}\beta(t,y)u(t,y) \d y
\\
&=&
-\int_{\Omega_B(t)}c(t,x_b) u(t,x) dx
-\int_{\Omega_B(t)}\partial_x c(t,x_b)(x - x_b ) u(t,x) dx
  +\sum_{i=B}^{J}\beta(t,x^i(t))m^i(t)
\\
&=&
-c(t,x_b)m^B(t) - \partial_x c(t,x_b)\pi^B(t)+\sum_{i=B}^{J}\beta(t,x^i(t))m^i(t),
\end{eqnarray*}
and furthermore, we expand $x c(x,t)$ around $x_b$
\begin{eqnarray*}
\mathrm{\frac{d}{dt}}\pi^B(t) &=&
\int_{\Omega_B(t)}(x_b - x) c(t,x)u(t,x) dx +
\int_{\Omega_B(t)} b(t,x)u(t,x) dx 
\\
&=&
x_b \int_{\Omega_B(t)}c(t, x_b )u(t,x) dx
+
x_b \int_{\Omega_B(t)} \partial_x c(t, x_b)(x - x_b) u(t,x) dx
\\
&&-\;
\int_{\Omega_B(t)} x_b c(t,x_b) u(t,x) dx
-
\int_{\Omega_B(t)} \Big(
c(t,x_b) + x_b\partial_x c(t, x_b)dx
\Big)(x - x_b)
+
\\
&& +\;
\int_{\Omega_B(t)} b(t,x_b)u(t,x) dx 
+
\int_{\Omega_B(t)} \partial_x b(t,x_b)(x - x_b)u(t,x) dx 
\\
&=&
b(t,x_b)m^i(t) + \partial_x b(t,x_b)\pi^B(t)
- c(t,x_b) \pi^B(t).
\end{eqnarray*}

\section{Proofs}
\begin{proofof}{Theorem \ref{existence_org}}
Proof of existence and uniqueness of solutions to
\eqref{EBT_org_system}~-~\eqref{EBT_org_BC} is not straightforward because of the specific definition of the dynamics of the boundary cohort \eqref{EBT_org_BC}, which implies that the right hand side is not Lipschitz continuous  on the whole $\R^{2(L+1)}$, where $L=J-B$. The term which causes difficulties is $\beta(t,\mu)(x^B(t))$ appearing in the last equation, since $x^B$ is given as a quotient $\pi^B/m^B$ and $m^B$ is not separated from zero. Therefore, we introduce a modification of the system whose right-hand side is Lipschitz continuous and then, prove that its solutions coincide with the solutions of the original system locally in time.

\begin{remark}
In the proof presented below, we assume that the dependence of the
model functions $b,c,\beta$ on the measure $\mu(t)$ is implicit, i.e.
\begin{equation}\label{additional}
\left\{
	\begin{array}{rcl}
\displaystyle
&&b(t, \mu) = \hat b \left(t, E_{b,\mu}\right),
c(t, \mu) = \hat c\left(t, E_{c,\mu}\right),
\beta(t, \mu) = \hat \beta\left(t, E_{\beta,\mu}\right),\\
&&E_{b,\mu} = \int_{\R_{\geq} } \gamma_b(y) \, d\mu(y),
E_{c,\mu} = \int_{\R_{\geq} } \gamma_c(y) \, d\mu(y),
E_{\beta,\mu} = \int_{\R_{\geq} } \gamma_{\beta}(y) \, d\mu(y)\\
&&\hat b, \hat c, \hat \beta \in C_{b}^{\alpha,1} \left(
 [0,T]\times\R_{\geq} ;\; W^{1,\infty}(\R_{\geq} ) \right),
\quad
\gamma_b, \gamma_c, \gamma_{\beta} \in W^{1,\infty} (\R_{\geq} ; \R_{\geq} ).
\end{array}
\right.
\end{equation}
The assumption about the implicit dependence of the model functions on the measure variable is not essential. The proof can be generalized by the use of standard, but quite technical, estimates.
\end{remark}

\textbf{Step 1: Problem on a restricted domain. }
Let us consider equation \eqref{EBT_org_system} with the following definition of the boundary cohort
\begin{eqnarray}\label{EBT_modified_BC}
\left\{
\begin{array}{rcl}
\displaystyle
x^B(t) &=&
\left\{
\begin{array}{lcl}
\min\{ \frac{\pi^B(t)}{m^B(t)}, C_K\}  + x_b,
&&
\displaystyle
\mathrm{if}\;\; m^B(t) > 0, \pi^B >0,
\\
 x_b,  && \mathrm{otherwise},
\end{array}
\right.
\\
\displaystyle
\mathrm{\frac{d}{dt}} \pi^B(t) &=& b(t, \mu)(x_b)m^B(t)
+ \partial_x b( t,\mu )(x_b)\pi^B(t)
\\
&& - c( t, \mu)(x_b)\pi^B(t),
\\[3mm]
\displaystyle
\mathrm{\frac{d}{dt}} m^B(t) &=&-c(t, \mu)(x_b)m^B(t)  - \partial_x c (t, \mu)(x_b)\pi^B(t) 
\\
&& + \sum_{i=B}^{J} \beta(t, \mu)(x^i(t))m^i(t),
	\end{array}
\right.
\end{eqnarray}
for some constant $C_K$.

The rand-hand side of system of the equations \eqref{EBT_org_system} and \eqref{EBT_modified_BC} is continuous, locally bounded and locally Lipschitz. Continuity and local boundedness are guaranteed by assumptions  \eqref{additional} which imply that the right-hand side can be estimated from above by
$$ C_1= C\left(\norma{(b,c,{\beta})}_{BC}, C_K,J-B\right)\sum_{i=B}^J \modulo{m^i}.$$

Lipschitz-continuity of the term $\beta(t,\mu)(x^B)$, which caused difficulties can now be proven as follows

\begin{eqnarray*}
&&|\beta(x^B)m^B- \beta(\tilde x^B)\tilde m^B| \leq \beta(x^B) |m^B-\tilde m^B| + \tilde m^B |\beta(x^B)- \beta(\tilde x^B)|\\
&&\leq \|\beta \|_{BC} \:\: |m^B-\tilde m^B| + \| \beta \|_{BC}\;\; \tilde m^B |x^B-\tilde x^B|
\\&&=  \|\beta \|_{BC}  \:\: |m^B-\tilde m^B| +  \|\beta \|_{BC} \;\; \left|\pi^B \left(1+ \frac{\tilde m^B-m^B}{m^B}\right) - \tilde \pi^B\right|\\
&&\leq \|\beta \|_{BC} \:\: |m^B-\tilde m^B| +  \|\beta \|_{BC} \left|\pi^B - \tilde \pi^B\right| +  \|\beta \|_{BC} \;\; |\tilde m^B-m^B| \: |\frac{\pi^B}{ m^B }|\\ 
&&\leq  \|\beta \|_{BC}  \:\: |m^B-\tilde m^B| +  \|\beta \|_{BC} \left|\pi^B - \tilde \pi^B\right| +  \|\beta \|_{BC} |\tilde m^B-m^B|  C_K
\end{eqnarray*}
on the set $\{(\pi^B,m^B):  \pi^B\leq C_K m^B\}$.
Because function $(\pi^B,m^B) \rightarrow x^B$ is Lipschitz-continuous on the whole domain, the above estimate extends to  $\{(\pi^B,m^B)\in \R^2 \}$.

Estimating $| E_{k, \mu} - E_{k, \tilde \mu}|$, we obtain term  $\gamma_k(x^B)m^B$ with $k\in \{b,c,\beta \} $, which can be handled exactly as the term $\beta(x^B)m^B$ above.

Taking into account all terms on the right-hand side yields its Lipschitz-continuity with a constant which can be estimated by
\begin{equation}
\label{c_lip}
C^* = C\left( 
\norma{(\gamma_b,\gamma_c,\gamma_{\beta}}_{W^{1,\infty}}, \norma{(\hat b, \hat c,\hat \beta)}_{BC}, C_K, J-B\right)\max\Big\{1, \sum_{i=B}^J \modulo{m^i}\Big\}.
\end{equation}
\noindent
Using Picard-Lindel\"of Theorem we conclude about local existence of solutions to problem \eqref{EBT_modified_BC}. The solution can be extended to  $[0,T]$ due to boundedness of the model coefficients $b, c, \beta$ and the fact that $\sum_{i=B}^{J}m^i(t)$ grows at most exponentially. 

\noindent
\textbf{Step 2: Well-posedness of the original problem}.
\\
We prove that there exists $T^*$ such that solutions of the original problem \eqref{EBT_org_system} - \eqref{EBT_org_BC} coincide with the solution of \eqref{EBT_org_system}, \eqref{EBT_modified_BC} for $t\in[0,T^*]$. Assuming existence of solutions of \eqref{EBT_org_system} - \eqref{EBT_org_BC}, we show that they are nonnegative and fulfill the restriction condition $\pi^B\leq C_K m^B$ locally in time. It yields local in time existence and uniqueness of solutions  of the original problem. 

For all times  such that $\sum_{i=B+1}^{J} \beta(t, \mu)(x^i(t))m^i(t)=0$, the solution $(\pi^B(t),m^B(t))$ remains equal to zero, so a positive $T^*$ trivially exists.  

Letting $\sum_{i=B+1}^{J} \beta(t, \mu)(x^i(t))m^i(t) > 0$, one needs to prove additional estimates to show that the (possibly) negative term $-\partial_x c(t,\mu)(x_b) \pi^B(t)$ never dominates the positive term $\sum_{i=B+1}^{J} \beta(t, \mu)(x^i(t))m^i(t)$. This can be done by bootstrapping argument. More precisely, one  integrates \eqref{EBT_org_BC} and applies the Gronwall's inequality. Substituting one of the resulting estimates into the other equation and iterating the procedure yields the following lemma.

\begin{lemma}\label{oszacowania}
There exist a time interval $[0,\tilde t]$ and constants $C_i>0$, which depend on $\tilde t $ and respective norms of model coefficients, such that the solutions of problem  \eqref{EBT_org_system} satisfy for $t\in[0,\tilde t ]$
\begin{eqnarray}\label{estimates_pi}
\modulo{\pi^B(t)}&\leq& C_1 t \int_0^{t} \sum_{i=B+1}^{J} \beta(\tau, \mu)(x^i(\tau))m^i(\tau) d\tau ,\\
m^B(t)&\geq& C_2 \int_0^{t} \sum_{i=B+1}^{J} \beta(\tau, \mu)(x^i(\tau))m^i(\tau) d\tau. 
\end{eqnarray}
\end{lemma}
From this we conclude about nonnegativity of $m^B$  and consequently also of $\pi^B$, for short enough time intervals.

To show that the trajectory remains in the cone $\{ (\pi^B,m^B): \pi^B \geq 0, m^B \geq 0, \pi^B\leq C_K m^B \}$, we follow \cite[Lemma 17]{EBT}.
The inequality $\pi^B\leq C_K m^B$ is fulfilled for $m^B(t_k) = 0$ and $\pi^B(t_k) = 0$.
Using Lemma \ref{oszacowania},
\begin{eqnarray}\label{limittk}
\label{der}
\frac{\mathrm{d}}{\mathrm{dt}}{x^B}|_{m^B=0}=\lim_{t\rightarrow t_k}\frac{\frac{\pi^B}{m^B}-0}{t-t_k}|_{m^B=0} \leq C_1/C_2.
\end{eqnarray}
For $t  \in (t_k,t_k+T^*)$ such that $m^B(t)>0$,  $x^B(t)$ is defined as a quotient of $m^B(t)/\pi^B(t)$ and thus, it is differentiable. We calculate
\begin{eqnarray}\label{xB}
\label{der}
\frac{\mathrm{d}}{\mathrm{dt}}\frac{\pi^B}{m^B}
&=&b(t,\mu)(x_b) + \partial_xb(t,\mu)(x_b) \frac{\pi^B(t)}{m^B(t)} 
 + \partial_x c(t,\mu)(x_b)\left(
 \frac{\pi^B(t)}{m^B(t)}
\right)^2 
\\
\nonumber
&&\quad
-  \frac{\pi^B(t)}{(m^B(t))^2} \sum_{i=B}^{J}\beta(t,\mu)(x^i(t)) m^i(t)
\\
\nonumber
&\leq&
b(t,\mu)(x_b) + \partial_xb(t,\mu)(x_b){x^B(t)} 
 + \partial_x c(t,\mu)(x_b)\left(
 {x^B(t)}
\right)^2.
\end{eqnarray}
It follows from \eqref{limittk} and \eqref{xB} that there exists $T^*$ such that
$
\frac{\mathrm{d}}{\mathrm{dt}}
x^B(t) \leq \bar C$  and consequently
$x^B(t) \leq \bar Ct$ for $t \in [t_k,t_k+T^*]$ and some constant $\bar C$.
In particular, this estimate implies that $x^B(\cdot)$ is continuous in $t = t_k$, which together with the continuity of the derivative on $(t_k,t_k+T^*)$, yields that $x^B(\cdot)$ is Lipschitz continuous on $[t_k,t_k+T^*]$. 
\end{proofof}


\begin{lemma}\label{LemmaFlatMetric}
Let $\mu = \sum_{i=1}^{L}m^i\delta_{x^i}$ and  $\tilde \mu = \sum_{i=1}^{L} \tilde m^i \delta_{\tilde x^i}$. Then,
\begin{eqnarray*}
\rho_F\left(\mu, \tilde \mu \right)
&\leq&
\max\Big\{
1, \sum_{i=1}^{L} \modulo{m^i}\Big\}
\sum_{i=1}^{L}( \modulo{x^i - \tilde x^i} +
\modulo{m^i - \tilde m^i}).
\end{eqnarray*}
\end{lemma}
\begin{proofof}{Lemma \ref{LemmaFlatMetric}}
\begin{eqnarray*}
\rho_F\left(\mu, \tilde \mu \right)
&\leq&
\rho_F\left( \sum_{i=1}^{L} m^i \delta_{x^{i}},  \sum_{i=1}^{L} m^i \delta_{\tilde x^{i}} \right)
+
\rho_F\left( \sum_{i=1}^{L} m^i \delta_{\tilde x^{i}},  \sum_{i=1}^{L} \tilde m^i \delta_{\tilde x^{i}}\right)\\
&\leq&W_1\left( \sum_{i=1}^{L} m^i \delta_{x^{i}},  \sum_{i=1}^{L} m^i \delta_{\tilde x^{i}} \right)
+
 \norma{ \sum_{i=1}^{L} m^i \delta_{\tilde x^{i}}-  \sum_{i=1}^{L} \tilde m^i \delta_{\tilde x^{i}}}_{TV}\\
 &\leq& 
\sum_{i=1}^{L}(
\modulo{m^i} \modulo{x^i - \tilde x^i} +
\modulo{m^i - \tilde m^i})
\\
\nonumber
&\leq&
\max\Big\{
1, \sum_{i=1}^{L} \modulo{m^i}\Big\}
\sum_{i=1}^{L}( \modulo{x^i - \tilde x^i} +
\modulo{m^i - \tilde m^i}),
\end{eqnarray*}
where $W_1$ is the $1$-Wasserstein distance and $\norma{\cdot}_{TV}$ is the total variation norm. 
\end{proofof}
\begin{proofof}{Theorem \ref{existence_org2}}
To prove existence and uniqueness of solutions we need to show that the right hand side is locally Lipschitz with respect to $(x^i, m^i)$. A proof of this claim can be conducted analogously as the corresponding part of the proof of Theorem \ref{existence_org}, therefore we do not repeat it here. The essential difference is that in instead of terms of the form
$
\modulo{E_{b,\mu} - E_{b,\tilde\mu}},
$
we obtain $\rho_F(\mu, \tilde \mu)$. Nevertheless, by Lemma \ref{LemmaFlatMetric} the right hand side is locally Lipschitz. Non-negativity of solutions is straightforward assuming that $b(t,\mu)(x_b) \geq 0$ and $\beta(t,\mu)(\cdot)$ is nonnegative. Exponential growth of mass can be proved analogously as in the proof of Theorem~\ref{existence_org}.
\end{proofof}

\section*{Acknowledgements}

{\small PG, AU and JJ were supported the International PhD Projects Programme of Foundation for Polish Science operated within the Innovative Economy Operational Programme 2007-2013 (PhD Programme: Mathematical Methods in Natural Sciences).  AU was additionally supported by the grant of National Science Centre no 2012/05/N/ST1/03132. AM-C was supported by the ERC Starting Grant No. 210680 and the Emmy Noether Programme of the German Research Council (DFG); and partially by National Science Centre (Narodowe Centrum Nauki).}

\end{document}